\documentclass[a4paper,3p,preprint]{elsarticle}
\usepackage{graphicx}
\usepackage{amsmath,amsfonts,latexsym,amssymb,amsbsy,amsthm}
\usepackage{natbib}
\usepackage{epstopdf}
\usepackage[makeroom]{cancel}
\usepackage{caption,subcaption}
\usepackage[multidot]{grffile}
\usepackage{bm}
\usepackage{booktabs}
\usepackage{framed}
\usepackage{mathtools}
\usepackage{hyperref,cleveref}

\def\O#1{{\cal O}\left(#1\right)}
\def\eps{\varepsilon}

\newcommand{\msp}{\phantom{-}}

\newcommand{\beq}{\begin{equation}}
\newcommand{\eeq}{\end{equation}}

\newcommand{\C}{\mathcal{C}}

\newcommand{\e}{\ensuremath{\mathrm{e}}}

\begin{document}

\begin{frontmatter}


	\title{Symplectic integrators for second-order linear non-autonomous equations}
	\journal{J. Comput. Appl. Maths.}

\author[lat]{Philipp Bader}\corref{cor1}\ead{p.bader@latrobe.edu.au}
\author[imm]{Sergio Blanes}\ead{serblaza@imm.upv.es}
\author[uji]{Fernando Casas}\ead{fernando.casas@uji.es}
\author[imm]{Nikita Kopylov}\ead{nikop1@upvnet.upv.es} 
\author[imm]{Enrique Ponsoda}\ead{eponsoda@imm.upv.es}
\address[lat]{Department of Mathematics and Statistics, La Trobe University, 3086 Bundoora VIC, Australia}
\address[imm]{Instituto de Matem\'{a}tica Multidisciplinar, 
Universitat Polit\`{e}cnica de Val\`{e}ncia. Spain}
\address[uji]{Universitat Jaume I, IMAC and Departament de Matem{\`a}tiques, 12071~Castell{\'o}n, Spain}
\cortext[cor1]{Corresponding author}



\begin{abstract}
   Two families of symplectic methods specially designed for second-order time-dependent linear systems are presented. Both are obtained from the Magnus expansion of the corresponding first-order equation, but otherwise they differ in significant aspects. The first family is addressed to problems with low to moderate dimension, whereas the second is more appropriate when the dimension is large, in particular when the system corresponds to a linear wave equation previously discretised in space. Several numerical experiments illustrate the main features of the new schemes.

\end{abstract}
\end{frontmatter}

	\section{Introduction}
	The problem we address in this paper is the~numerical integration of the~second order time-dependent linear equation
	\begin{equation}\label{hill}
	x''(t) + M(t) x(t) = 0 \, , \qquad x(0)=x_0,\quad x'(0)=x_0',
	\end{equation}
	where $ t \in \mathbb{R} $, $ x(t) \in \mathbb{C}^{r} $ and primes denote time-derivatives. It
	has many practical applications, for example, 
	in periodically variable systems such as qua\-dru\-po\-le mass filter and~qua\-dru\-po\-le devices, 
	\cite{drewsen00hlp,paul90etf}, 
	microelectromechanical systems \cite{turner98fpr}, 
	in Bose--Einstein condensates,
	spatially linear electric fields, dynamic buckling of structures, electrons in crystal lattices, waves in periodic media, etc. (see \cite{bader16sif,magnus66he,major05cpt,mclachlan65taa} and~references therein).
	In these cases the~matrix $M(t)$ is usually time periodic with period $T$ ((\ref{hill}) is then an example of a matrix Hill equation) and of moderately large size $r$. Parametric resonances can occur so that it is very important to~know the~stability regions in terms of the~parameters of the~system. 
	
	Another important example of this type is the linear time-dependent wave equation 
	\begin{equation}\label{eq.se}
	\partial_t^2 u(x,t) = f(x,t)\partial_x^2 u(x,t)+ g(x,t)u(x,t),\quad x\in\mathbb{R},\ t\geq0,
	\end{equation}
	equipped with initial conditions $u(x,0)=u_0(x)$, and~$u_t(x,0)=u'_0(x)$. Once discretized in space in a bounded region, \cref{eq.se} also leads to  \cref{hill}.
	
	When the~dimension $r$ of \cref{hill} is relatively small, the~numerical computation of its fundamental matrix solution for one period (usually repeatedly many times for different values of the~parameters of the~system) is feasible and allows one to~analyse the stability of the~configuration. In this case numerical schemes that involve matrix--matrix products can be used. However, when \cref{hill} results from a semidiscretized PDE like \cref{eq.se}
	then $r\gg 1$, and so only numerical schemes involving matrix--vector products are suitable.
	
	Taking these considerations into account, in 
	this work we present two classes of numerical schemes: one that involves matrix--matrix products showing a~high performance for the~computation of the~fundamental matrix solution for one period, and another class involving matrix--vector products, addressed for the numerical treatment of PDEs.

	Equation (\ref{hill}) can be written as a~first order system  by introducing new variables $q=x, p=x'$ as
	%
	\begin{equation}\label{firstordersystem}
	z'(t) = A(t) z(t) \, , \qquad \mbox{ with } \qquad 
	A(t) = \begin{pmatrix}
	0 & I \\
	- M(t) & 0
	\end{pmatrix} ,
	\end{equation}
	and $z=(q,p)^T$, $ \, z_0 =(q_{0},p_{0})^T\in\mathbb{C}^{2r}$. 
	The solution of \cref{firstordersystem} evolves through a~linear transformation (evolution operator or fundamental matrix solution) given by $z(t)=\Phi(t,0)z(0)$.
	In the~usual case in which $M(t) $ is a~real and symmetric $ r \times r $ matrix valued function ($M^T=M$) then $\Phi(t,0)$ is a~symplectic transformation.
	The eigenvalues of $\Phi(t,0)$ occur in reciprocal pairs, say $\{\lambda,\lambda^*, 1/\lambda, 1/\lambda^*\}$, where $\lambda^*$ denotes the~complex conjugate of $ \lambda $. As a~result, for stable systems, all of the~eigenvalues must lie on the~unit circle (see for example \cite{dragt15lmf}). This is a~very important property that is not preserved in general by standard numerical integrators. Thus, one can be forced to~use very small time steps to~avoid the undesirable numerical instabilities or asymptotic stabilities arising from this fact, with the resulting degradation in the efficiency of the algorithms.

	The solution of the non-autonomous \cref{firstordersystem} cannot be written in general in closed form. 
	Nevertheless, we can always formally write it as a~single exponential by using the~Magnus expansion \cite{blanes09tme,magnus54ote}, 
	\[ \begin{pmatrix} q(t_n+h)\\p(t_n+h)\\ \end{pmatrix} = 
	e^{\Omega(t_n,h)}\begin{pmatrix} q(t_n)\\p(t_n)\\ \end{pmatrix}, \]
	where $\Omega(t_n,\tau)$ 
	satisfies a (highly nonlinear) differential equation. Although different approximations can be found in the literature 
	(see e.g. \cite{blanes09tme} and references therein), they involve the computation of nested commutators of $A(t)$, evaluated at different times,
	and so for many problems may be computationally expensive\footnote{We must remark however that in some cases, like for the~linear Schr\"odinger equation with time-dependent potential, the~evaluation of the~commutators can be efficiently carried using appropriate approach \cite{bader16emf}.} 
	
	For this reason, here we consider two different strategies starting from the formal solution provided by the Magnus expansion. The first leads us
	to optimised methods specifically addressed to problems defined by (\ref{hill}) when matrix--matrix products are feasible; to these we call
	\textit{Magnus-decomposition methods}. The second strategy allows us to
	build optimised methods when only matrix--vector products are suitable. These are referred to as \textit{Magnus--splitting methods}. 
	

	In \textit{Magnus-decomposition methods}, in particular,  
	we approximate $\exp\Omega(t_n,h)$ by a~composition of simpler symplectic maps that for the problem at hand are considerably faster to~compute.
	As an illustration, a 4th-order method within this family contains  only one costly matrix exponential and this exponential is cheaply approximated by a symplectic approximation of order $q$ leading to an scheme given by the composition:
	\begin{equation}\label{eq:phi4q1}
	\Upsilon_1^{[4,q]} = \begin{pmatrix} I & 0 \\ D_{3}^{[4,q]} & I \end{pmatrix}
	\begin{pmatrix} 0 & D_{2}^{[4,q]} \\ I & 0 \end{pmatrix}
	\begin{pmatrix} I & 0 \\ D_{1}^{[4,q]} & I \end{pmatrix},
	\end{equation}
	where 
	$D_i^{[4,q]}, \ i=1,2,3$,  
are certain matrices to be defined later on.
	
	
	On the other hand, 
	in \textit{Magnus--splitting methods}, we consider different decompositions that are suitable for low-oscillatory and high-dimensional problems, where only matrix--vector products are reasonable in the~numerical scheme. Methods of order $p$ within this family are of the form  
	\begin{equation}\label{eq:RKN2}
	\Psi_m^{[p]} = 
	\begin{pmatrix} I & h a_{m+1} I \\ 0 & I \end{pmatrix} 
	\begin{pmatrix} I & 0 \\ h C_m & I \end{pmatrix} 
	\begin{pmatrix} I & h a_m I \\ 0 & I \end{pmatrix} 
	\ \cdots
	\begin{pmatrix} I & 0 \\ h C_1 & I \end{pmatrix} 
	\begin{pmatrix} I & a_1 I \\ 0 & I \end{pmatrix},
	\end{equation}
	where the~coefficients $a_i,b_i$ are obtained numerically by solving a~set of non-linear polynomial equations and the matrices $C_i$ are
	linear combinations of $M(t)$ evaluated at certain quadrature points.
	The~product $\Psi_m^{[p]} z_0$ is done with only $m$ products of matrices $C_i$ on a~vector of dimension $r$.
	
	Since all the~schemes presented in this work are based on the Magnus expansion, they become explicit symplectic integrators when $M$ is real and symmetric and can be considered as geometric integrators \cite{blanes16aci,hairer06gni,leimkuhler05shd,sanzserna94nhp}, showing a favourable behaviour on long-time simulations.  
	
	\section{Magnus-based methods}
	

	\subsection{Magnus expansion}
	
	As stated before, the Magnus expansion \cite{magnus54ote} expresses the~solution of ~\cref{firstordersystem} on the time interval [$t_n$,\,$t_n+h$] in the~form of a~single exponential of an infinite series
	\begin{equation}\label{Magnus}
	\Phi(t_n,h) = \exp \Omega (t_n,h), \qquad 
	\Omega (t_n,h) = \sum_{k=1}^{\infty} \Omega_{k} (t_n,h) ,
	\end{equation}
	whose first terms are given by
	\begin{equation}\label{omegaks}
	\Omega_{1} (t_n,h) = 
	\int\limits_{{t_n}}^{t_n+h}\!A(t_{1}) \, d\tau_{1} \, , \quad
	\Omega_{2} (t_n,h) = \frac{1}{2} \int\limits_{{t_n}}^{t_n+h} \int\limits_{t_n}^{\tau_{1}} \left[ A(\tau_{1}), A(\tau_{2}) \right] \, d\tau_{2} \, d\tau_{1}, 
	\end{equation}
	where $[\alpha,\beta] = \alpha\beta - \beta\alpha$ is the~matrix commutator.  When matrix $A(t)$ is given by \eqref{firstordersystem}, the exponent
	$\Omega$ and as well as any truncation of the series at order $ p $, $\Omega^{[p]}$, belong to~the symplectic Lie algebra, and
	thus symplecticity is automatically preserved, 
	
	Approximations of $\Omega$ in terms of $A(t)$ evaluated at the~nodes of some quadrature rule can be obtained as follows.
	First we consider the~polynomial $\widetilde A(t)$ of degree $s-1$ in $t$ that interpolates $A(t)$ on $[t_n,t_n+h]$ at the~points $t_n+c_jh, \ j=1,\ldots,s$, where $c_i$ are the~nodes of the~Gauss--Legendre quadrature rule of order $2s$.
	The perturbed problem reads
	\begin{equation} \label{eq:ivpInt} 
	\frac{d \widetilde z(t)}{d t} = \widetilde A(t) \, \widetilde z(t), \qquad \widetilde z(t_n)=z(t_n), \qquad
	t \in [t_n,t_n+h],
	\end{equation}
	where $z(t_n)$ is the~exact solution of (\ref{firstordersystem}) at $t_n$.
	From a~direct application of the~Alekseev--Gr\"obner lemma \cite{zanna99car} (see also \cite{iserles99ots,moan98eao,munthekaas99cia}), we have that
	\[
	\widetilde z(t_n+h)-z(t_n+h) = {\cal O}(h^{2s+1}).
	\]
	Letting $t=t_n+\frac{h}2+\sigma$, we write the~interpolation polynomial as
	\begin{equation} \label{tilde_A}
	\widetilde A(t) = \sum_{i=1}^s {\cal L}_i\Big(\frac{t-t_n}{h}\Big)A_i
	=\frac1{h}\sum_{i=1}^s 
	\left(\frac{\sigma}{h}\right)^{i-1} \alpha_i,\qquad\quad 
	\sigma\in\left[-\frac{h}2,\frac{h}2\right],
	\end{equation}
	with $A_i=A(t_n+c_i h)=\widetilde A(t_n+c_i h)$ and~the usual Lagrange polynomials ${\cal L}_i(t)$.
	Notice that for our problem (\ref{firstordersystem}) we have
	\begin{equation} \label{alpha_i}
	\alpha_{i+1} = \left. h^{i+1}\frac{1}{i!} \, 
	\frac{d^{i} \widetilde A(t_{1/2}+\sigma)}{d \sigma^{i}} \right|_{\sigma=0},
	\qquad i=0,\ldots,s-1,
	\end{equation}
	where 
	\[
	\alpha_1=\begin{pmatrix} 0 & I \\ -\mu_1 & 0 \end{pmatrix},\qquad\quad
	\alpha_j=\begin{pmatrix} 0 & 0 \\ -\mu_j & 0 \end{pmatrix},\qquad j>1,
	\]
	and
	\begin{equation} \label{mu_i}
	\mu_{i+1} = \left. h^{i+1}\frac{1}{i!} \, 
	\frac{d^{i} \widetilde M(t_{1/2}+\sigma)}{d \sigma^{i}} \right|_{\sigma=0}\!,
	\qquad i=0,\ldots,s-1.
	\end{equation}
	In particular, for sixth-order methods ($s=3$) we have
	\begin{equation} \label{alpha_iA_i}
	\alpha_1 = h A_2, \qquad \alpha_2 = \frac{\sqrt{15}h}{3} (A_3 - A_1), \qquad
	\alpha_3 = \frac{10h}{3} (A_3 - 2 A_2 + A_1),
	\end{equation}
	where $A_i=A(t_n+c_ih), \ c_1=\frac{5-\sqrt{15}}{10}, \ c_2=\frac12, \ c_3=\frac{5+\sqrt{15}}{10}$. One can readily check that $\alpha_1={\cal O}(h), \alpha_2={\cal O}(h^2),\alpha_3={\cal O}(h^3)$.
	
	The~integrals in the~Magnus expansion for $\widetilde A(t)$ can be computed from \eqref{eq:ivpInt} analytically. 
	This results in an~approximation of~$z(t)$ expressed in terms of $\alpha_1,\: \alpha_2,\:\alpha_3$ up to~order $2s$.
	Specifically,  a sixth-order approximation $\Omega^{[6]} = \Omega + \mathcal{O} (h^{7})$ results from
	\begin{equation} \label{sixth}
	\Omega^{[6]} =\alpha_{1} + \frac{1}{12} \alpha_{3} - \frac{1}{12} [12] +
	\frac{1}{240} [23] +
	\frac{1}{360} [113] - \frac{1}{240} [212] + \frac{1}{720} [1112],
	\end{equation}
	where $[ij\ldots kl]$ represents the~nested commutator $\lbrack
	\alpha_{i},[\alpha_{j}, [\ldots ,[\alpha_{k},\alpha_{l}]\ldots]]]$.

	At this point, we can proceed in different ways and~we present two different strategies to~build new methods. In both cases, the respective
	methods (\ref{eq:phi4q1}) and (\ref{eq:RKN2}) will be compositions of exponentials of elements of the Lie algebra generated by	
	$\alpha_1,\: \alpha_2,\:\alpha_3$ whose coefficients will be determined in such a way that the compositions coincide with (\ref{sixth}) up to the desired
	order (six in this case).
	
	\subsection{Magnus-decomposition integrators}

	Firstly, we examine the~structure of the~Lie algebra generated by the~$\alpha_i$. We immediately notice that, $[ij]=0$ for $i,j>1$.
	Furthermore, 
	\[
	[212]=\begin{pmatrix} 0 & 0 \\ 2 \mu_2^2 & 0 \end{pmatrix}.
	\]
	
	We distinguish the~following types of exponentials that may appear as elements of the~Lie algebra generated by $\alpha_i, \ i=1,2,3$:
	\begin{equation}
	\label{E1E2E3}
	E_1=\exp\begin{pmatrix}D & B\\ C & -D^T\end{pmatrix}, \qquad
	E_2=\exp\begin{pmatrix}0 & I\\ C & 0\end{pmatrix}, \qquad
	E_3=\exp\begin{pmatrix}0 & 0\\ C & 0\end{pmatrix}.
	\end{equation}
	Clearly, we want to~avoid the~computation of the~full matrix exponential $E_1$ and~instead focus on types $E_2,\:E_3$. The~latter is nilpotent and~its exponential comes virtually for free.

	In \cite{bader16sif}, it was obtained methods of order 4, 6 and 8 for the~Hill equation, denoted by $\Upsilon_k^{[p]}$ (a 
	$ p^{th} $-order method containing $ k $ exponentials of type $E_2$ above). In particular	 
	\begin{equation}\label{eq:phi41}
	\Upsilon_1^{[4]} =\begin{pmatrix} I & 0 \\ hC_{2}^{[4]} & I \end{pmatrix} \,\, 
	\exp \left(h \begin{pmatrix} 0 & I \\ D_{1}^{[4]} & 0 \end{pmatrix}\right)  \,\,
	\begin{pmatrix} I & 0 \\ hC_{1}^{[4]} & I \end{pmatrix}
	\end{equation}
	\begin{equation}\label{eq:phi62}
	\Upsilon_2^{[6]} =\begin{pmatrix} I & 0 \\ hC_{2}^{[6]} & I \end{pmatrix} \,\,
	\exp \left(\frac{h}{2} \begin{pmatrix} 0 & I \\ D_{2}^{[6]} & 0 \end{pmatrix}\right)  \,\,
	\exp \left(\frac{h}{2} \begin{pmatrix} 0 & I \\ D_{1}^{[6]} & 0 \end{pmatrix}\right)  \,\,
	\begin{pmatrix} I & 0 \\ hC_{1}^{[6]} & I \end{pmatrix},
	\end{equation}
	where $ C_{i}^{[p]}$~are linear combinations of $M(t)$ evaluated at a~set of quadrature points of order $p$ or higher and $ D_{i}^{[p]}$~are linear combinations of $M(t)$ that additionally contain one product of such linear combinations. 
	Due to the~exponentials, these schemes are specially appropriate when the~solution is oscillatory or stiff.
	
	Since the~averaged matrices  $D_{i}^{[p]}$ can be considered as constant matrices in each time subinterval, 
	the following result is useful when computing the corresponding exponentials for a real-valued matrix $C=-M$ \cite[sec. 11.3.3]{Bernstein2009}:
	\begin{equation} \label{eq:solaut}
	\Phi(h) = \exp\left( h \begin{pmatrix}
	0 & I \\ 
	C & 0
	\end{pmatrix} \right) = 
	\begin{pmatrix}
	\cosh h\sqrt{C}     & \sqrt{C}^{-1}\sinh h\sqrt{C} \\ 
	\sqrt{C} \sinh h\sqrt{C} & \cosh h\sqrt{C}       
	\end{pmatrix}. 
	\end{equation}
	
	Although there exist efficient methods to~compute matrix trigonometric functions \cite{almohy15naf,alonso16eaa} as well as matrix exponentials, such as Pad\'e approximants, 
	Krylov/Lanczos methods, Chebyshev method and~others, for this particular case a~more efficient procedure is obtained by decomposing 
	the~exponential into a~product of simple matrices. When $C$ is symmetric, the resulting 
	approximations preserve the~symplectic structure by construction. Specifically, if $h\rho(\sqrt{C})<\pi$, where $\rho(\sqrt{C})$ 
	is the~spectral radius of $\sqrt{C}$, then \cite{bader14stp}
	\begin{equation} 
	\Phi(h) = \begin{pmatrix}
	I & 0 \\ 
	R & I
	\end{pmatrix} 
	\begin{pmatrix}
	I & Q \\ 
	0 & I 
	\end{pmatrix} 
	\begin{pmatrix}
	I & 0 \\ 
	R & I 
	\end{pmatrix} ,
	\label{eq:factoring}
	\end{equation}
	where
	\begin{align}\label{eq:factoring:sintan}
	\begin{split}
		 Q(C)&=\frac{\sinh h \sqrt{C}}{\sqrt{C}} = hI+\frac{C h^3}{6}+\frac{C^2 h^5}{120}+\frac{C^3 h^7}{5040}+\frac{C^4 h^9}{362880}+\frac{C^5 h^{11}}{39916800}+\frac{C^6 h^{13}}{6227020800}+\O{h^{15}} \\
	R(C)&=\sqrt{C} \tanh\left(\frac{h\sqrt{C}}{2}\right) = \frac{C h}{2}-\frac{C^2 h^3}{24}+\frac{C^3 h^5}{240}-\frac{17 C^4 h^7}{40320}+\frac{31 C^5 h^9}{725760}-\frac{691 C^6h^{11}}{159667200}+ \O{h^{13}} 
	\end{split}		
	\end{align}
	The~truncated series expansions of $Q$ and $R$ up to order $q+2$ and $ q $ in $h$, denoted by $Q^{[q+2]}$ and $R^{[q]}$, respectively,
	can be simultaneously computed with only $k=\left\lfloor \frac{q-1}2\right\rfloor$ products (e.g., $Q^{[6]},\: P^{[4]}$ can be computed with only one product).  We denote by $\Phi^{[q]}$ a $ q^{th} $-order approximation to~$\Phi(h)$ obtained by replacing $ Q $ and~$ R $ by $Q^{[q+2]}$ and $R^{[q]}$. 
	Notice that if $C$ is a~symmetric matrix then $Q^{[q+2]},\:R^{[q]}$ are also symmetric matrices and, by construction, $\Phi^{[q]}$ is a~symplectic matrix $\forall q$. 
	Taking into account all these considerations, we substitute \cref{eq:factoring,eq:factoring:sintan} into $\Upsilon_k^{[p]}$ given by 
	(\ref{eq:phi41}) and (\ref{eq:phi62})  and combine commuting matrices to get finally $\Upsilon_{1}^{[4,q]}, \Upsilon_{2}^{[6,q]}$ given in Table~\ref{tab:methods41_62}.
	The algorithm proceeds by multiplying $\Upsilon_{k}^{[p,q]}$  by the~result from the previous computational step. As the 
	last matrix in a step commutes with the~first one in the~following step, some products can be saved, hence for an~even $ q\geq6 $ the 
	computational cost of 
	$ \Upsilon^{[4,q]}_{1} $ and $ \Upsilon^{[6,q]}_{2} $ are $ (1+q/2)\C $ and $ (7+q)\C$, respectively. 
	The schemes and the relevant parameters are collected in \Cref{tab:methods41_62}.
	
	\begin{table}
		\caption{One- and two-exponential 4th- and 6th-order symplectic method, respectively, using the~sixth-order Gauss--Legendre quadrature rule.}
		\label{tab:methods41_62}
		\begin{framed}
			\[
			c_1=\frac{1}{2}-\frac{\sqrt{15}}{10}, \quad c_2=\frac12, \quad c_3=\frac{1}{2}+\frac{\sqrt{15}}{10}.
			\]
			\[
			M_1=M(t_n+c_1h), \qquad M_2=M(t_n+c_2h), \qquad M_3=M(t_n+c_3h), 
			\]
			\[
			K=M_1-M_3, \qquad L=-M_1+2M_2-M_3, \qquad F=h^2 K^2.
			\]
			\begin{center}
				$  Q_i^{[p,q+2]}$ and $R_i^{[p,q]} $ are obtained by applying the expansions \cref{eq:factoring:sintan} to $ D_i^{[p]} \!.$
			\end{center}
			\begin{framed}
				\[
				\begin{array}{ll}
				C_1^{[4]}=\frac{\sqrt{15}}{36}K + \frac{5}{36}L \qquad
				&  
				D_1^{[4]}=-M_2 \\
				C_2^{[4]}=-\frac{\sqrt{15}}{36}K + \frac{5}{36}L
				&  
				\end{array}
				\]
				\begin{equation*}
				\Upsilon_{1}^{[4,q]}=
				\begin{pmatrix}
				I           & 0 \\ 
				hC_2^{[4]}+R_1^{[4,q]} & I 
				\end{pmatrix} 
				\begin{pmatrix}
				I & Q_1^{[4,q+2]} \\ 
				0 & I       
				\end{pmatrix} 
				\begin{pmatrix}
				I           & 0 \\ 
				hC_1^{[4]}+R_1^{[4,q]} & I 
				\end{pmatrix}.
				\end{equation*}
			\end{framed}
			\begin{framed}
				\[
				\begin{array}{ll}
				C_1^{[6]}=-\frac{\sqrt{15}}{180}K + \frac1{18}L+\frac1{12960}F \qquad
				&  
				D_1^{[6]}=-M_2 -\frac{4}{3\sqrt{15}}K + \frac1{6}L \\
				C_2^{[6]}=+\frac{\sqrt{15}}{180}K + \frac1{18}L+\frac1{12960}F
				&  
				D_2^{[6]}=-M_2 +\frac{4}{3\sqrt{15}}K + \frac1{6}L
				\end{array}
				\]
				\begin{equation*}
				\Upsilon_{2}^{[6,q]}=
				\begin{pmatrix}
				I & 0 \\ 
				hC_2^{[6]}+R_2^{[6,q]} & I
				\end{pmatrix} 
				\begin{pmatrix}
				I & Q_2^{[6,q+2]} \\ 
				0 & I
				\end{pmatrix} 
				\begin{pmatrix}
				I & 0 \\ 
				R_2^{[6,q]}+R_1^{[6,q]} & I
				\end{pmatrix} 
				\begin{pmatrix}
				I & Q_1^{[6,q+2]} \\ 
				0 & I
				\end{pmatrix} 
				\begin{pmatrix}
				I & 0 \\ 
				hC_1^{[6]}+R_1^{[6,q]} & I
				\end{pmatrix}.
				\end{equation*}
			\end{framed}
		\end{framed}
	\end{table}

	\subsection{Magnus--splitting integrators}
	For deriving the second class of schemes considered in this work, we first split the~matrix $A(t)$ of \cref{firstordersystem} as 
	\begin{eqnarray} \label{eq.AB} 
	A(t)=B(t) + D \qquad \mbox{ with } \qquad 
	B(t) = \left( \begin{array}{cc}
	0 & 0 \\ -M(t) & 0 \end{array}
	\right), \qquad D = \left( \begin{array}{cc}
	0 & I \\ 0 & 0 \end{array}
	\right),
	\end{eqnarray}
	and denote
	\[
	\delta_1 = hD, \qquad\quad
	\displaystyle \beta_i = \frac{h^i}{(i-1)!} \,
	\frac{d^{i-1} \widetilde B(s)}{d s^{i-1}} \Big|_{s=t+\frac{h}{2}}, \ i \geq1
	\]
	where $\widetilde B(s)$ is the~interpolating polynomial of~$B(s)$ in the~interval $[t_n, t_n + h]$,
	and $\alpha_1=\delta_1+\beta_1$ $\alpha_i=\beta_i, \ i>1$. 
	It is easy to~check that $[\beta_i,\beta_j]=0$ and
	$[\delta_1,\delta_1,\delta_1,\beta_i]=[\beta_i,\beta_j,\beta_k,\delta_1]=0$
	for any value of $i,j,k$. As a consequence, the formal solution \cref{sixth} simplifies to
	\begin{eqnarray} \label{A-constant}
	\Omega^{[6]} &=& \delta_1 + \beta_1 +
	\frac{1}{12} \beta_3 + \frac{1}{12}
	[\beta_2,\delta_1]+ \frac{1}{360} \big( -[\delta_1,\beta_3,\delta_1] + [\beta_1,\delta_1,\beta_3] \big) \nonumber\\
	& & - \frac{1}{240} [\beta_2,\delta_1,\beta_2] + \frac{1}{720} \big([\delta_1,\beta_1,\delta_1,\beta_2]
	- [\beta_1,\delta_1,\beta_2,\delta_1] \nonumber
	\big).
	\end{eqnarray}
	
	It is then clear that a composition of type (\ref{eq:RKN2}) can be recovered by considering the following product of exponentials:
	\begin{equation}  \label{comp.1}
	\Psi_m^{[6]}  =  \prod_{i=1}^{m+1} \exp\left(\sum_{j=1}^3 y_{i,j}\beta_j 
	\right) \, \exp\left(a_i\delta_1\right)
	\end{equation}
$y_{m+1,j}=0, \ j=1,2,3$.	Specifically, taking to account (\ref{eq.AB}), we have
	\begin{equation}  \label{comp.2}
	\Psi_m^{[6]}  =  \prod_{i=1}^{m+1} 
	\begin{pmatrix} I & 0 \\ \displaystyle \sum_{j=1}^3b_{i,j}h M_j & I \end{pmatrix}
	\begin{pmatrix} I & a_ihI \\ 0 & I \end{pmatrix}
	\end{equation}
	where $M_j=M(t_n+c_jh), \ j=1,2,3$ and $b_{m+1,j}=0$, so that, in practice, \eqref{comp.2} corresponds to~a~$m$-stage composition.
	All methods are symplectic when applied to~Hamiltonian systems and, moreover, the coefficients are chosen so that time symmetric is preserved.
	
	To~obtain particular methods we extend the analysis carried out in \cite{blanes07smf,blanes12smi}, where several 6th-order schemes were derived for the~more general problem
	\begin{equation}   \label{eq.ho}
	q'=M(t)p, \qquad p'=N(t)q.
	\end{equation}
	Although all of them can be applied to~the present problem \eqref{firstordersystem} simply by letting $M(t)=I$, new methods have also been obtained by taking into consideration the
	simpler structure that this system possesses. Specifically, we have taken  $a_i,\:y_{i,1}$ in \eqref{comp.1} as the coefficients of an optimised 11-stage 6th-order method 
	designed in \cite{blanes12smi} for \eqref{eq.ho} when $M$ and $N$ are constant. In this way the~commutators involving only $\delta_1$ and $ \beta_1$  (e.g. $[\delta_1,\beta_1,\delta_1]$,
	$[\beta_1,\delta_1,\beta_1]$, etc.)  vanish up to~order six, whereas the~higher-order contributions are minimised by considering more stages than strictly necessary to~solve the~order conditions; this extra cost is compensated by a~much improved accuracy and stability in an autonomous case. Next, we look for new coefficients $y_{i,2},\:y_{i,3}$, which now have to satisfy a much 
	reduced set of order conditions. The 11-stage composition has five  coefficients $y_{i,2}$ and six coefficients $y_{i,3}$ to solve seven equations. Taking into account the structure of the equations, this leaves one of the  $y_{i,2}$ and three of the  $y_{i,3}$ as free parameters, that are chosen in order to minimize the objective function $\sum_i (|y_{i,2}|^2+|y_{i,3}|^2)$.
	
	Once the~coefficients $y_{i,j}$ are chosen, the~matrices $\beta_i, \ i=1,2,3$ are replaced by the~corresponding linear combinations of $M(t_n+c_ih), \ i=1,2,3$ evaluated at the~quadrature rule of order 6, so that one ends up with a composition of the form (\ref{comp.2}). Specifically, 
	the~following 11-stage 6th-order method is obtained:
	\begin{equation}\label{eq:psi2_11}
	\Psi_{11}^{[6]} = 	\begin{pmatrix} I & h a_{12} I \\ 0 & I \end{pmatrix} 
	\begin{pmatrix} I & 0 \\ hC_{11} & I \end{pmatrix} 
	\begin{pmatrix} I & ha_{11}I \\ 0 & I \end{pmatrix}
	\ \cdots
	\begin{pmatrix} I & 0 \\ hC_{1} & I \end{pmatrix} 
	\begin{pmatrix} I & ha_1I \\ 0 & I \end{pmatrix},
	\end{equation}
	where
	\[
	C_i = -(b_{i,1} M_1 + b_{i,2} M_2 + b_{i,3} M_3), \qquad i=1,\ldots, 11,
	\]
	$M_j=M(t_n+c_jh)$ and $c_i$ are the~nodes of the~6th-order Gauss--Legendre quadrature rule, i.e., 
	\[
	c_1=\frac{5-\sqrt{15}}{10}, \quad c_2=\frac12, \quad c_3=\frac{5+\sqrt{15}}{10}.
	\]   
	The corresponding coefficients $a_i, b_{i,j}$ are:
	\[
	\begin{array}{lll}
	a_1= 0.04648745479086313 &  
	a_2=-0.06069167116564293 &  
	a_3= 0.21846652646340681 \\
	a_4= 0.16805357948309270 &  
	a_5= 0.31439236417035348 &  
	a_6=-0.18670825374207319 
	\end{array}
	\]
	
	{
		\begin{equation} \label{11stage_bs}
		\big(b_{i,j}\big) = \begin{pmatrix}
		\msp 0.152309756970167&  \msp  0.078927889445323&  -0.046907162912825\\
		\msp 0.006406269275594 &-0.091413523927685&   \msp 0.043950351354379\\
		\msp 0.086778862327312  & \msp 0.051027214890409 & -0.004050397550970\\
		\msp 0.066634120201024   &\msp 0.148499347182669 & -0.011368920251338\\
		-0.020231991304321   &	\msp 0.030206484536889 & -0.021734660147529\\
		\msp 0.025991549816284   &\msp 0.009949620189233 &  \msp 0.025991549816284\\
%
		
		\end{pmatrix}
		\end{equation}
	}
	verifying the time-symmetry  condition $ a_{13-i} = a_{i}, \quad i=1,\ldots,6$, and 
	\[
	b_{6+i,j} = b_{6-i,4-j}, \qquad i=1,\ldots,5, \ \ j=1,2,3.
	\]
	
	\section{Numerical experiment}
	In the~following section we study and~demonstrate the~performance of the~new methods with respect to~well established explicit and implicit standard Runge--Kutta (RK) and explicit symplectic Runge--Kutta--Nystr\"om (RKN) methods from the~literature. Six types of methods are considered:
	\begin{itemize}
		\item $ \Upsilon_k^{[p]} $ from \cite{bader16eni} and \cite{bader16sif}: symplectic $ p^{th} $-order methods requiring computing of $ k $ matrix exponentials.
		\item $ \Upsilon_{k}^{[p,q]} $: new sets of methods obtained from $ \Upsilon^{[p]}_k $ by decomposing matrix exponentials and taking Taylor series expansion up to~$ \O{h^q} $.
		\item $ \Psi_{11}^{[6]} $: the new 11-stage 6th-order Magnus--splitting method (\cref{eq:psi2_11,11stage_bs}).
		\item $ \widehat{\Psi}_{11}^{[6]} $: an 11-stage 6th-order Magnus--splitting method for non-autonomous systems from \cite{blanes07smf,blanes12smi}.
		\item $ \mathrm{RK}_{k}^{[p]}$: an~explicit $ k $-stage $ p^{th} $-order method that uses the~Radau quadrature rule of order six for the~time-dependent part.
		\item $ \mathrm{RKGL}_{s}^{[p]}$: implicit $ k $-stage $ p^{th} $-order symplectic Gauss--Legendre methods.
		\item $ \mathrm{RKNb}_{k}^{[p]}$: $ k $-stage $ p^{th} $-order explicit symplectic methods from \cite{blanes02psp}. $ \mathrm{RKNb}_{11}^{[6]} $ has been selected instead of the more effective $ \mathrm{RKNa}_{14}^{[6]} $ (ibid.) because the former has the same number of stages as $ \Psi_{11}^{[6]} $.
	\end{itemize}
	
	
	The cost of each method is estimated for two different problem types. The~first one is when a~numerical method acts on the~fundamental matrix $ \Phi $. Let $r\times r =\dim M(t)$.  Then, the~cost is expressed in the~number of  matrix\,--\,matrix products $ \C $, required to~propagate for one time step $ h $. Evaluations of $ M(t) $, scalar\,--\,matrix multiplications and linear combinations of matrices are not included to~the cost.
	
	The second case is when the same method is used to~integrate a~system whose state is represented by a~vector $ (v,w)^T, \dim v=\dim w =r $. Similarly, the~cost is expressed in the~number of matrix\,--\,vector products $ \mathcal{V} $. The~methods' costs are summarized in Table~\ref{tab:cost}. The parameter $ \varrho^{[p]} $ in the~implicit methods refers to~the~number of iterations per step for a $p^{th}$-order method. Typically, $ \varrho^{[p]}=4...7$ to~attain convergence and preservation of symplecticity to~high accuracy. For the numerical experiments in this paper, we assume them to be $ \varrho^{[4]}=4 $ for the 4th-order method and $ \varrho^{[6]}=6 $ for the 6th-order one.

	\begin{table}[h]
		\centering
		\caption{Cost of the~methods in terms of matrix--matrix products, $ \C $, and matrix--vector products, $ \mathcal{V} $. }\label{tab:cost}
		\begin{tabular}[t]{lcc}
			\toprule
			Method           & $ \C $       & $ \mathcal{V} $  \\ 
			\cmidrule{1-3}
			$ \Upsilon^{[4]}_1 $    & $ 17\frac{1}{3} $ & --         \\ 
			$ \Upsilon^{[4,6]}_{1} $  & 4         & --         \\ 
			$ \Upsilon^{[4,8]}_{1} $  & 5         & --         \\ 

			$\mathrm{RK}^{[4]}_{4} $  & 8         & 4         \\ 
			$\mathrm{RKGL}^{[4]}_{2}$ & $4\times \varrho^{[4]}$ & $2\times \varrho^{[4]}$ \\ 
			$\mathrm{RKNb}^{[4]}_{6} $ & 12         & 6         \\ 
			\bottomrule
		\end{tabular}
	\qquad
		\begin{tabular}[t]{lcc}
			\toprule
			Method          & $ \C $       & $ \mathcal{V} $  \\ 
			\cmidrule{1-3}
			$ \Upsilon^{[6]}_2 $   & $ 32\frac{2}{3} $ & --         \\ 
			$ \Upsilon^{[6,8]}_{2} $ & 15         & --         \\ 
			$ \Upsilon^{[6,12]}_{2} $ & 19         & --         \\ 
			$ \Psi_{11}^{[6]},\, \widehat{\Psi}_{11}^{[6]} $    & 22         & 11         \\ 

			$\mathrm{RK}^{[6]}_{7} $  & 14         & 7         \\ 
			$\mathrm{RKGL}^{[6]}_{3} $ & $6\times \varrho^{[6]}$ & $3\times \varrho^{[6]}$ \\ 
			$\mathrm{RKNb}^{[6]}_{11} $ & 22         & 11         \\ 
			\bottomrule
		\end{tabular}
	\end{table}

	 However, some of these methods are not optimised for problems where matrix--matrix are exceedingly costly, and they are not used in these problems. Moreover, the~schemes $\Upsilon_{k}^{[p]}$ acting on a~vector could be carried out using propagators like Krylov methods to~compute the~action of the~exponential of a~matrix on a~vector, but this is not considered in this work. 
	
	The reference solutions are obtained numerically using sufficiently small time steps.
	
	\subsection{Mathieu equation}
	The first performance test is executed employing the~Mathieu equation
	\begin{equation}\label{ex:mathieu}
		x''(t)+(\omega^2+\eps\cos 2t)x(t)=0,
	\end{equation}
	written as a~first-order system (\ref{firstordersystem}).
		
	At first, we family-wise compare the~performance of the~new methods $ \Upsilon^{[p,q]}_{k} $. 	
	Then, we integrate for $ t\in [0,\pi ] $ with the~identity matrix as the~initial condition and~then measure the~$ L_1 $-norm of the~error of the~fundamental matrix at the~final time. This procedure is repeated for different time steps and different choices of $k,p,q$. In Table~\ref{tab:best_q_of_omega} we show the values of $q$ that provided the best performances for the choices $\eps=\{0.1,\:1\}$ and $\omega=5^j, \ j=-3,\ldots,3$. We observe that the~choices $ q=6,\,8 $ are generally better in the 4th-order family, and $q=8,\,12$ are better among the 6th-order methods. These selected schemes will be considered in the~following numerical examples.
	 \begin{table}
		\caption{The orders of decomposition of the two best performing methods. The better one comes first.}
		\label{tab:best_q_of_omega}	
		\centering
		\begin{tabular}{cccccccc}
			\toprule
			 $\omega$  & 1/125 & 1/25 & 1/5  & 1   & 5   & 25  & 125  \\ 
			\hline
			\rule[-1ex]{0pt}{4ex}
			& \multicolumn{7}{c}{$\varepsilon=1$} \\
			\rule[-1ex]{0pt}{3ex} $p$ of $\Upsilon_{1}^{[4,p]}$ & 6, 10 & 6, 10 & 6, 10 & 6, 10 & 12, 10 & 8, 12 & 8, 12 \\ 
			\rule[-1ex]{0pt}{2.5ex} $q$ of $\Upsilon_{2}^{[6,q]}$ & 10, 8 & 10, 8 & 10, 8 & 10, 8 & 12, 8 & 12, 8 & 12, 8 \\ 
			\rule[-1ex]{0pt}{4ex}
			& \multicolumn{7}{c}{$\varepsilon=1/10$} \\
			\rule[-1ex]{0pt}{2.5ex} $p$ of $\Upsilon_{1}^{[4,p]}$ & 10, 6 & 10, 6 & 10, 6 & 10, 12 & 12, 8 & 8, 12 & 8, 12 \\ 
			\rule[-1ex]{0pt}{2.5ex} $q$ of $\Upsilon_{2}^{[6,q]}$ & 8, 10 & 8, 10 & 8, 10 & 8, 10 & 12, 8 & 12, 8 & 12, 8 \\
			\bottomrule
		\end{tabular}
	\end{table}

\begin{figure}[!h]	
		\begin{subfigure}[b]{0.5\textwidth}	
			\centering
			\includegraphics[width=1\linewidth]{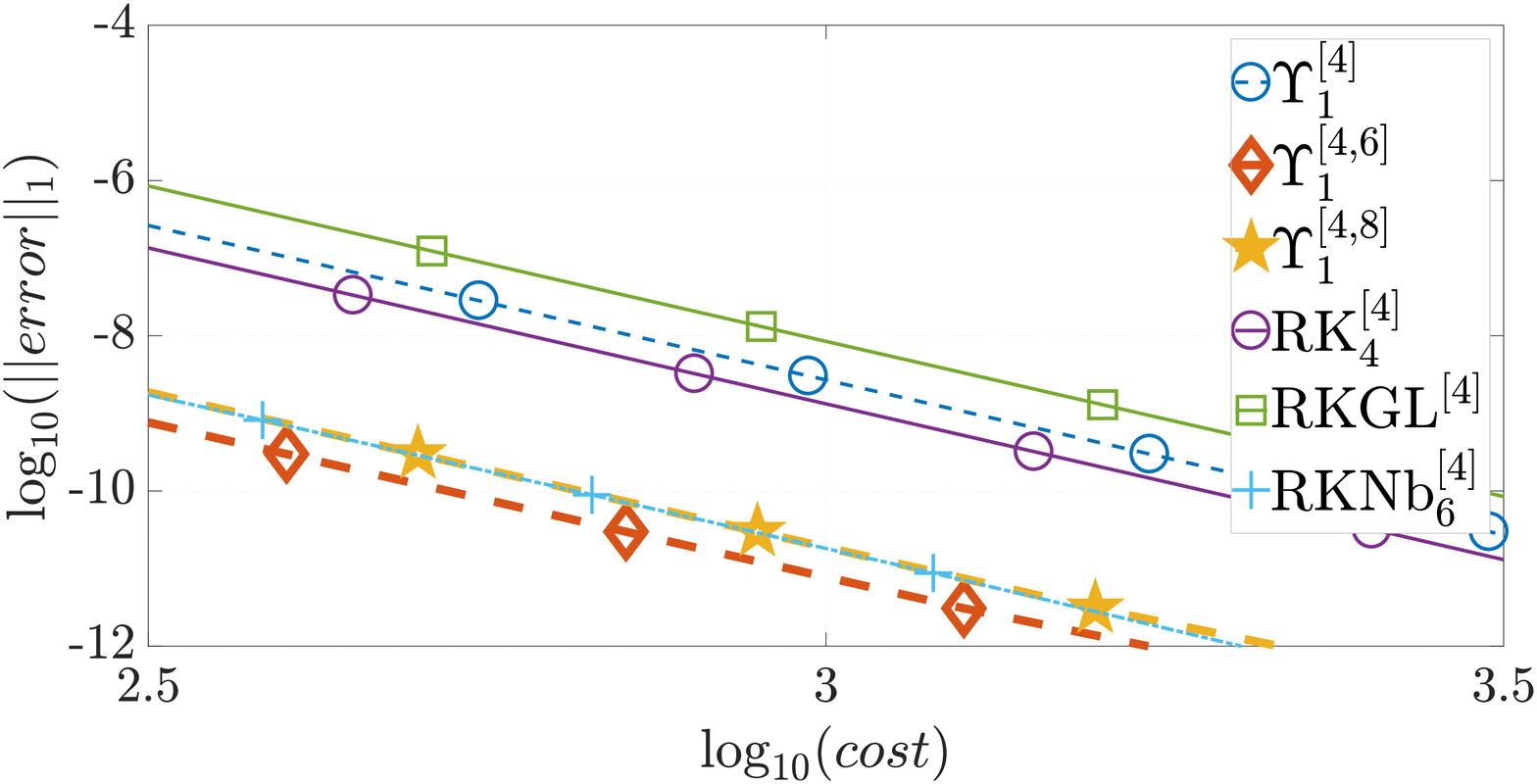}
			\caption{$ \omega=1/5,\,\varepsilon=1/10 $}
		\end{subfigure}
		\begin{subfigure}[b]{0.5\textwidth}	
			\centering
			\includegraphics[width=1\linewidth]{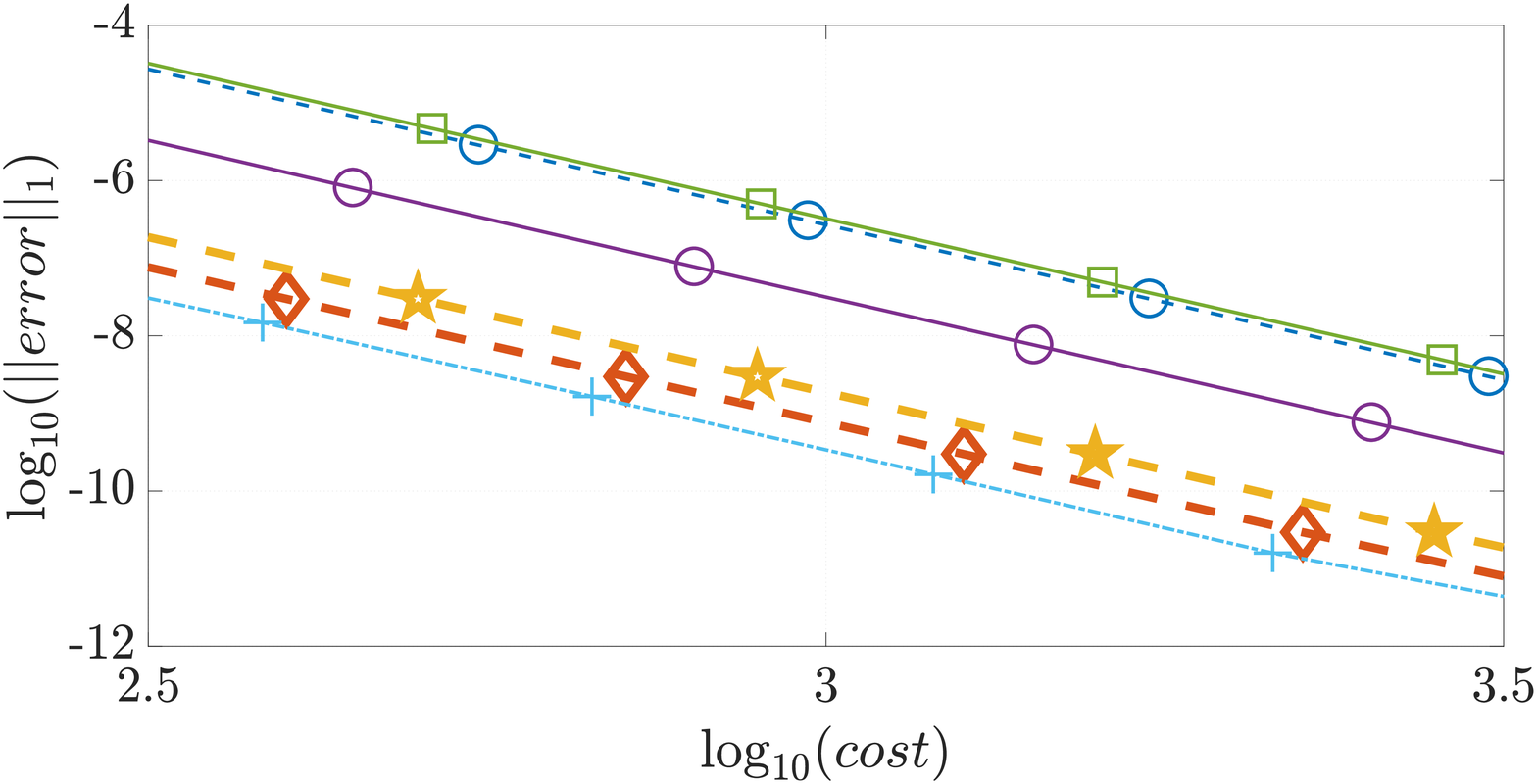}
			\caption{$ \omega=1/5,\,\varepsilon=1 $}
		\end{subfigure}
		\begin{subfigure}[b]{0.5\textwidth}
			\centering
			\includegraphics[width=1\linewidth]{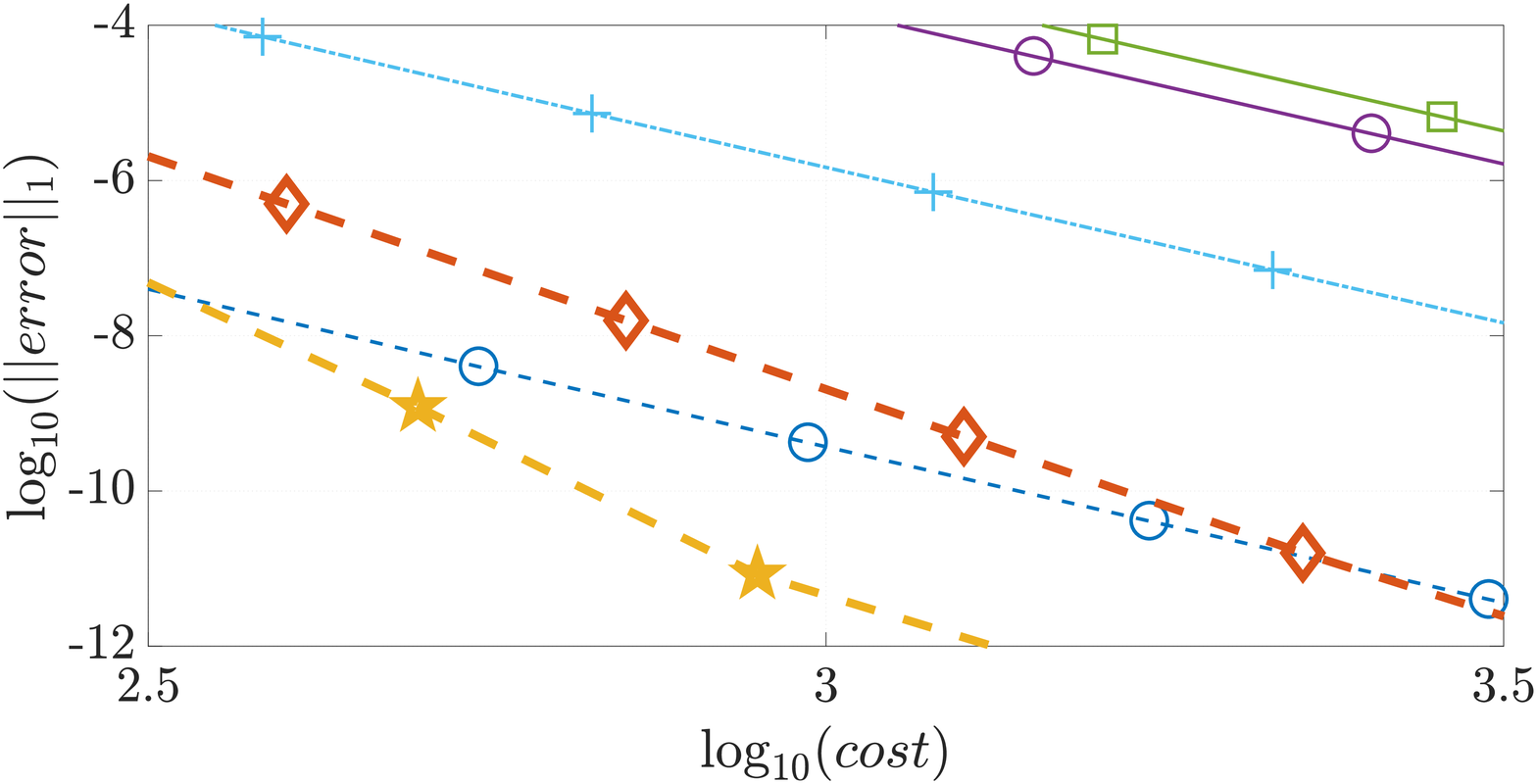}
			\caption{$ \omega=5,\,\varepsilon=1/10 $}
		\end{subfigure}
		\begin{subfigure}[b]{0.5\textwidth}
			\centering
			\includegraphics[width=1\linewidth]{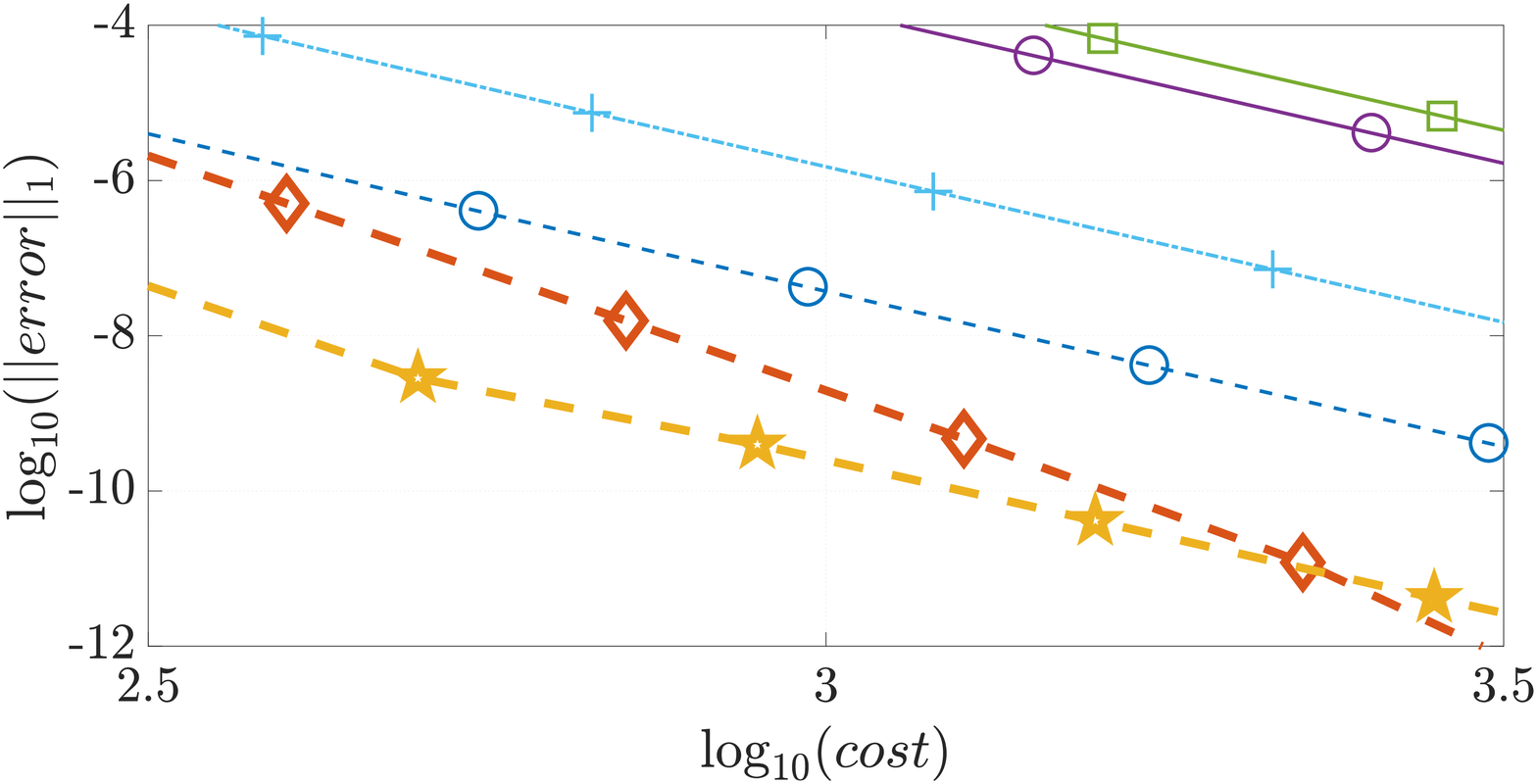}
			\caption{$ \omega=5,\,\varepsilon=1 $}
		\end{subfigure}
		\caption{The performance of the 4th-order methods for the Mathieu \cref{ex:mathieu} on a $ \log $--$ \log $ scale; $ cost = \mathcal{C}\times steps $.}
		\label{fig:mth_4th}
	\end{figure}

\begin{figure}[!h]	
		\begin{subfigure}[b]{0.5\textwidth}	
			\centering
			\includegraphics[width=1\linewidth]{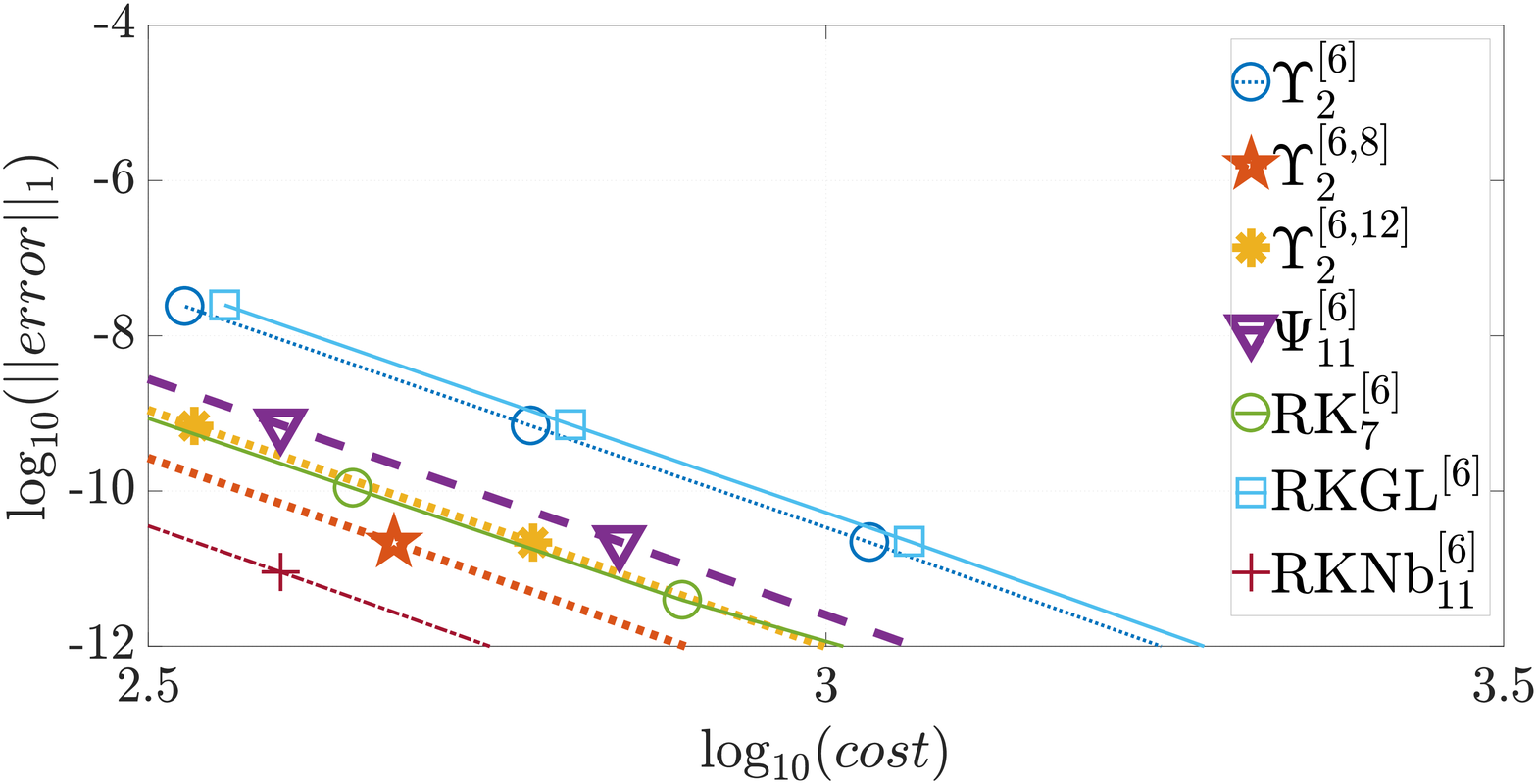}
			\caption{$ \omega=1/5,\,\varepsilon=1/10 $}
		\end{subfigure}
		\begin{subfigure}[b]{0.5\textwidth}	
			\centering
			\includegraphics[width=1\linewidth]{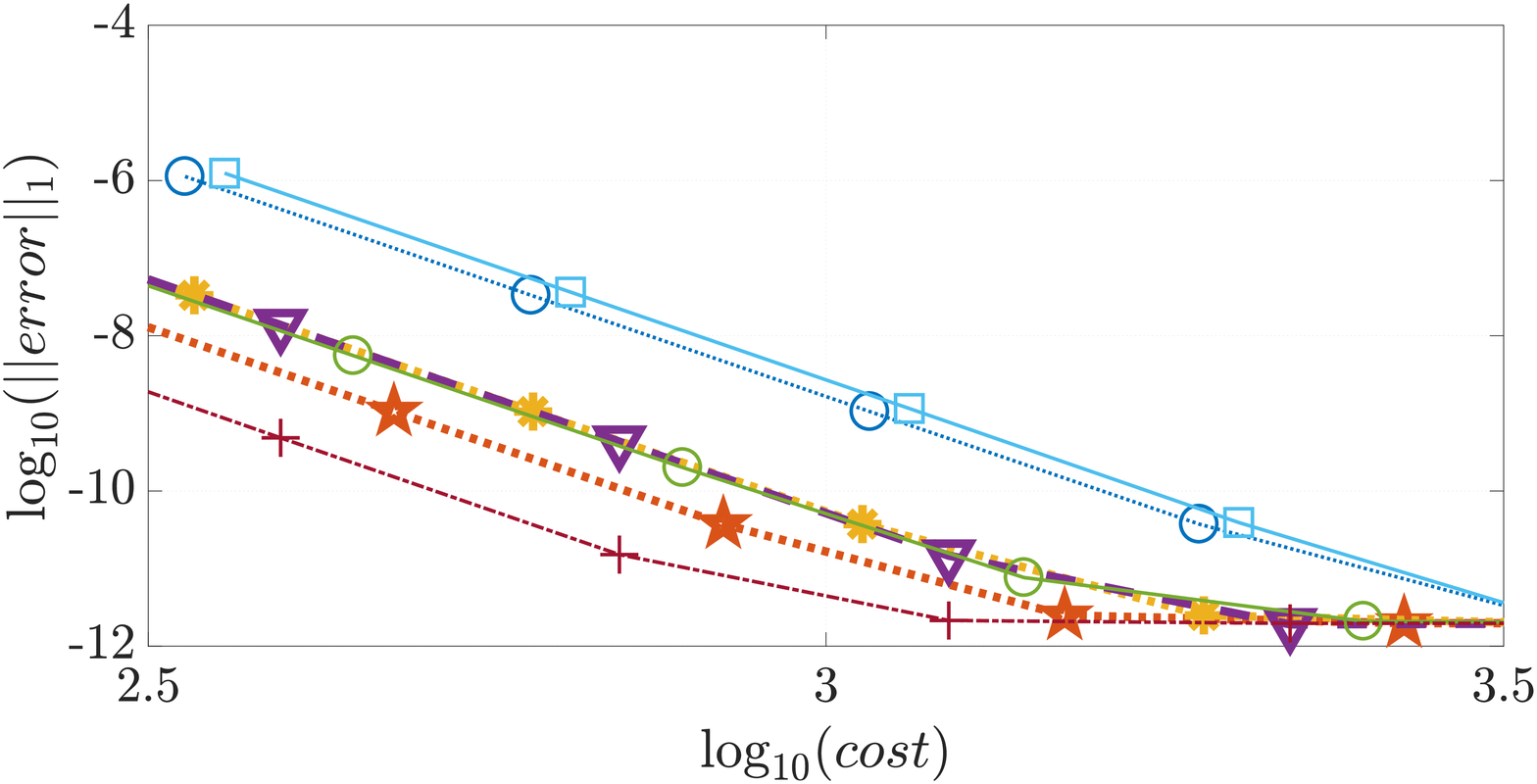}
			\caption{$ \omega=1/5,\,\varepsilon=1 $}
		\end{subfigure}
		\begin{subfigure}[b]{0.5\textwidth}
			\centering
			\includegraphics[width=1\linewidth]{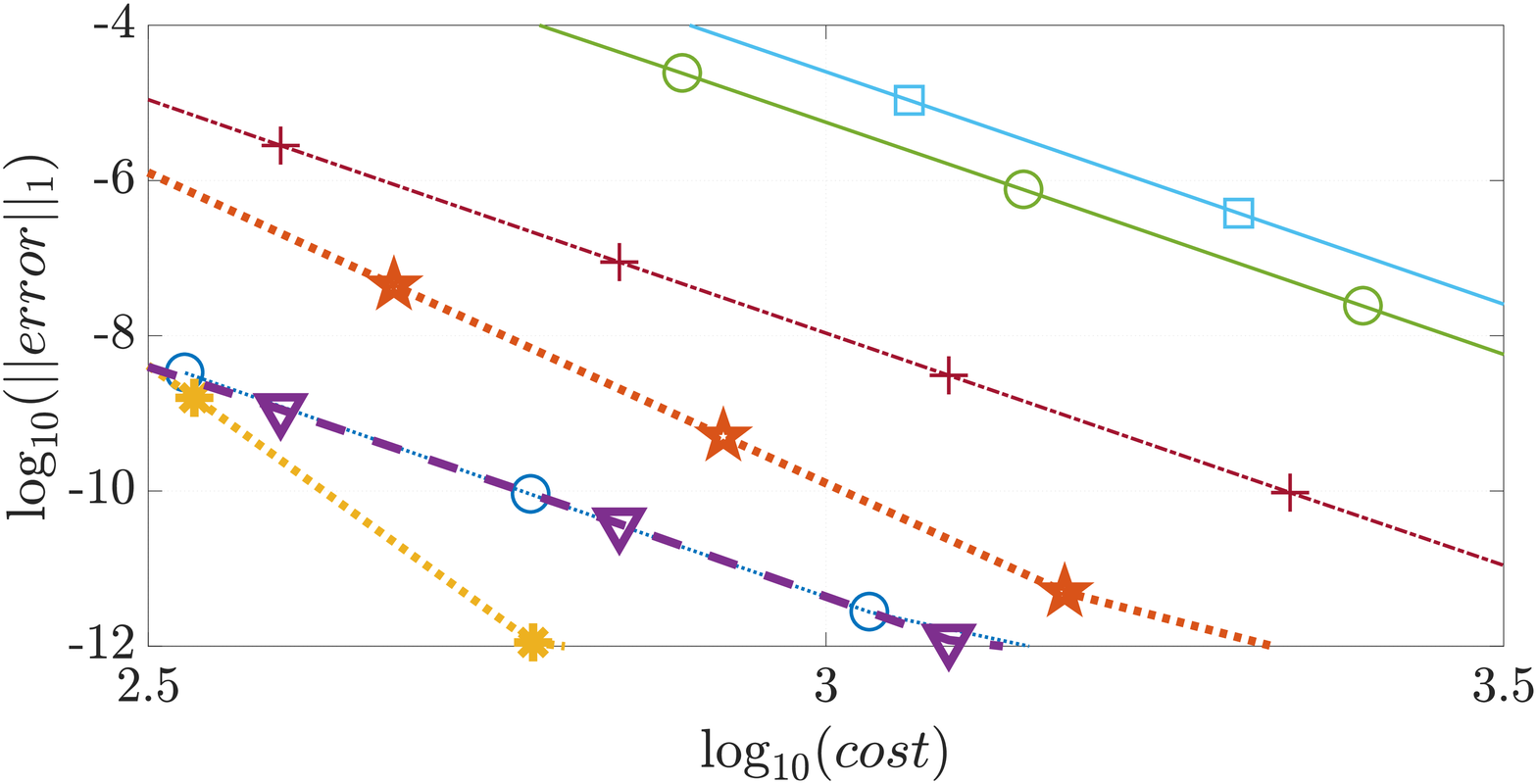}
			\caption{$ \omega=5,\,\varepsilon=1/10 $}
		\end{subfigure}
		\begin{subfigure}[b]{0.5\textwidth}
			\centering
			\includegraphics[width=1\linewidth]{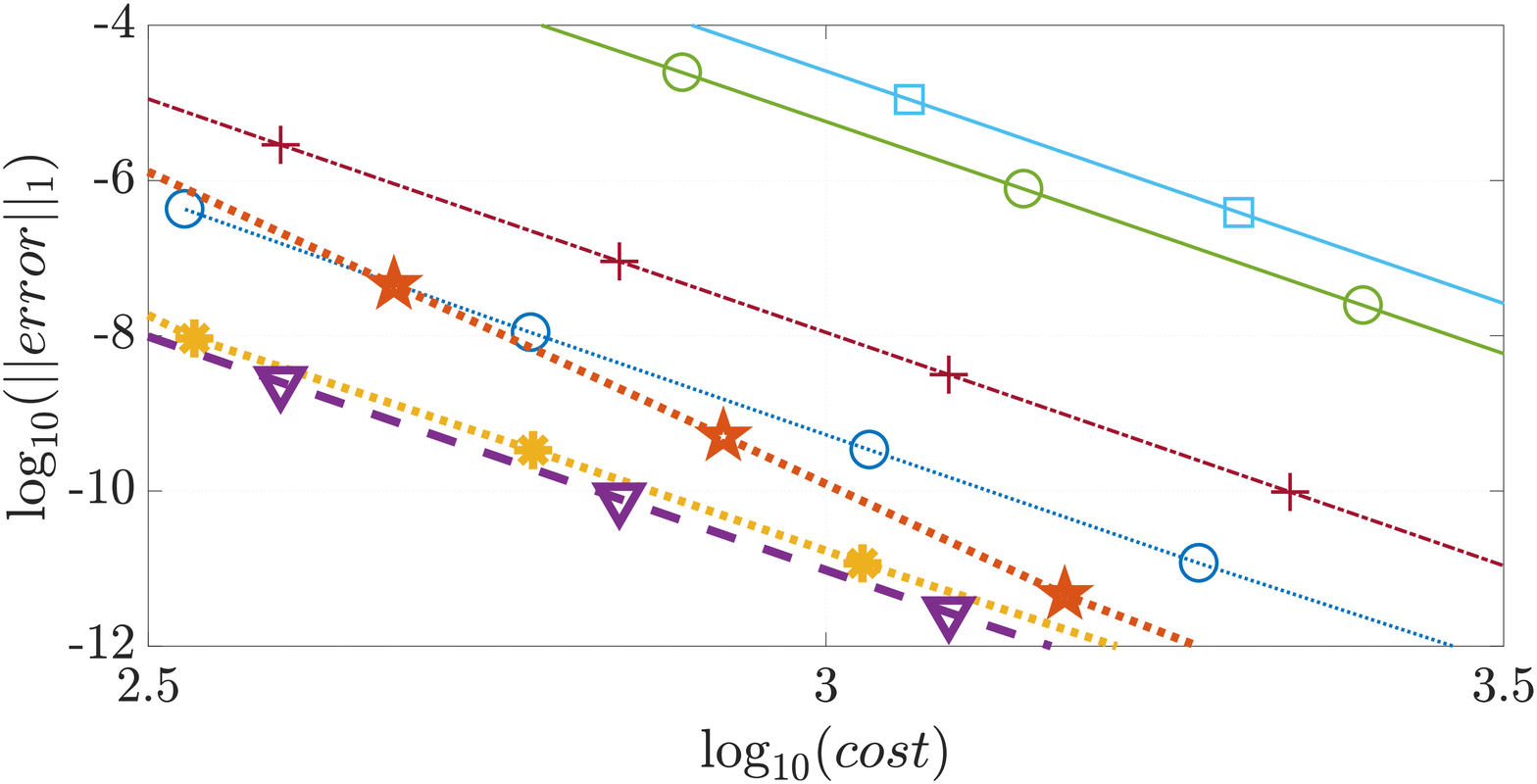}
			\caption{$ \omega=5,\,\varepsilon=1 $}
		\end{subfigure}
		\caption{The performance of the 6th-order methods for the Mathieu \cref{ex:mathieu} on a $ \log $--$ \log $ scale; $ cost = \mathcal{C}\times steps $.}
		\label{fig:mth_6th}
	\end{figure}

	We select now $ \eps\in \{0.1,\: 1\}$ and~$ \omega\in\{1/5,\: 5\} $ for plotting the~$ L_1 $-norm of the~error of the~fundamental matrix at the~final time versus the~computation cost in units of $ \C $ for all methods previously considered.
	\Cref{fig:mth_4th,fig:mth_6th} show methods' performance.
	We observe that the~standard explicit and implicit RK methods perform considerably worse while the~new methods show the~best performances for oscillatory cases.	
	
\begin{figure}[!h]		
		\begin{subfigure}[b]{0.5\textwidth}	
			\centering
			\includegraphics[width=1\linewidth]{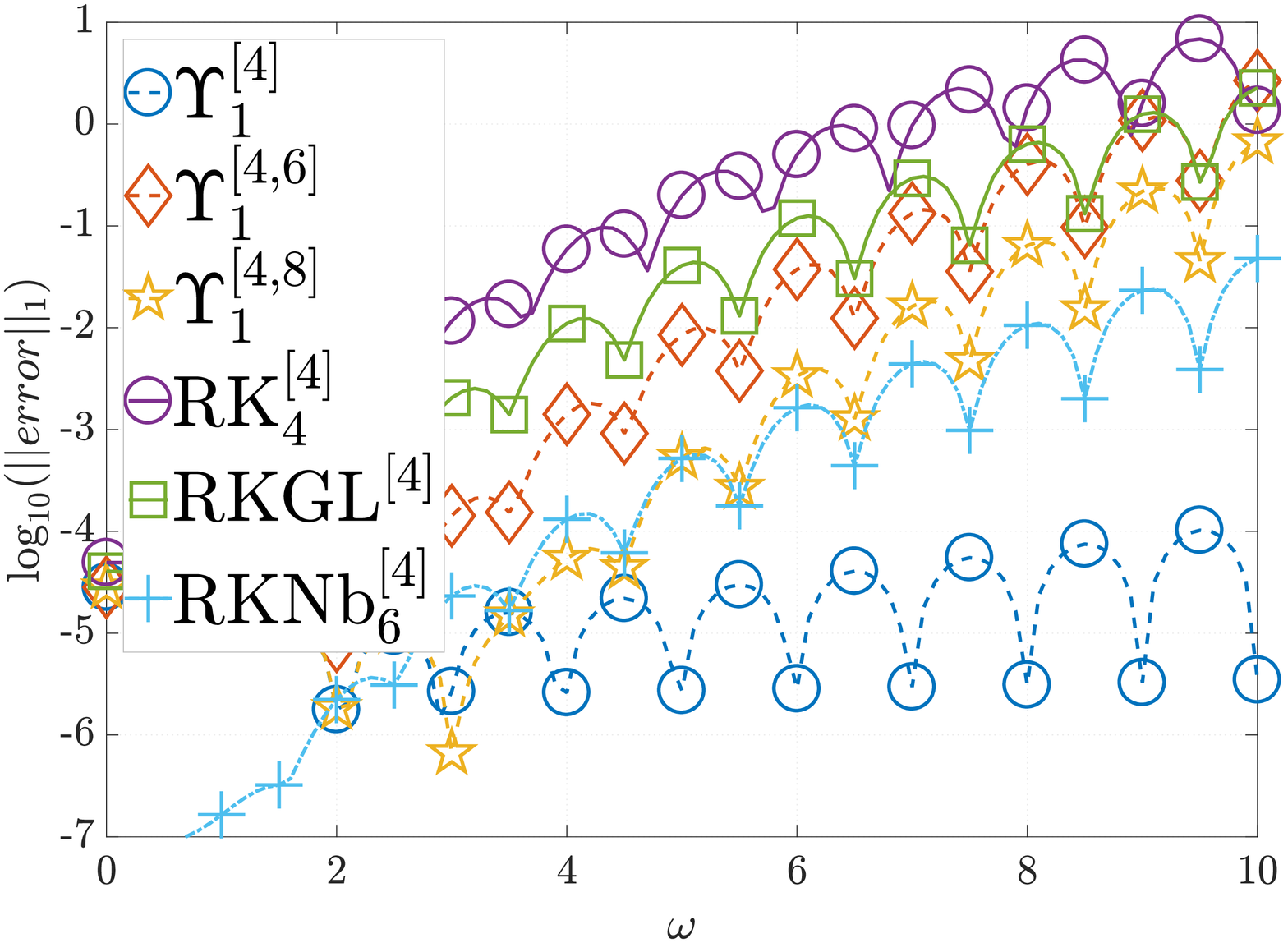}
			\caption{4th-order methods}
		\end{subfigure}
		\begin{subfigure}[b]{0.5\textwidth}	
			\centering
			\includegraphics[width=1\linewidth]{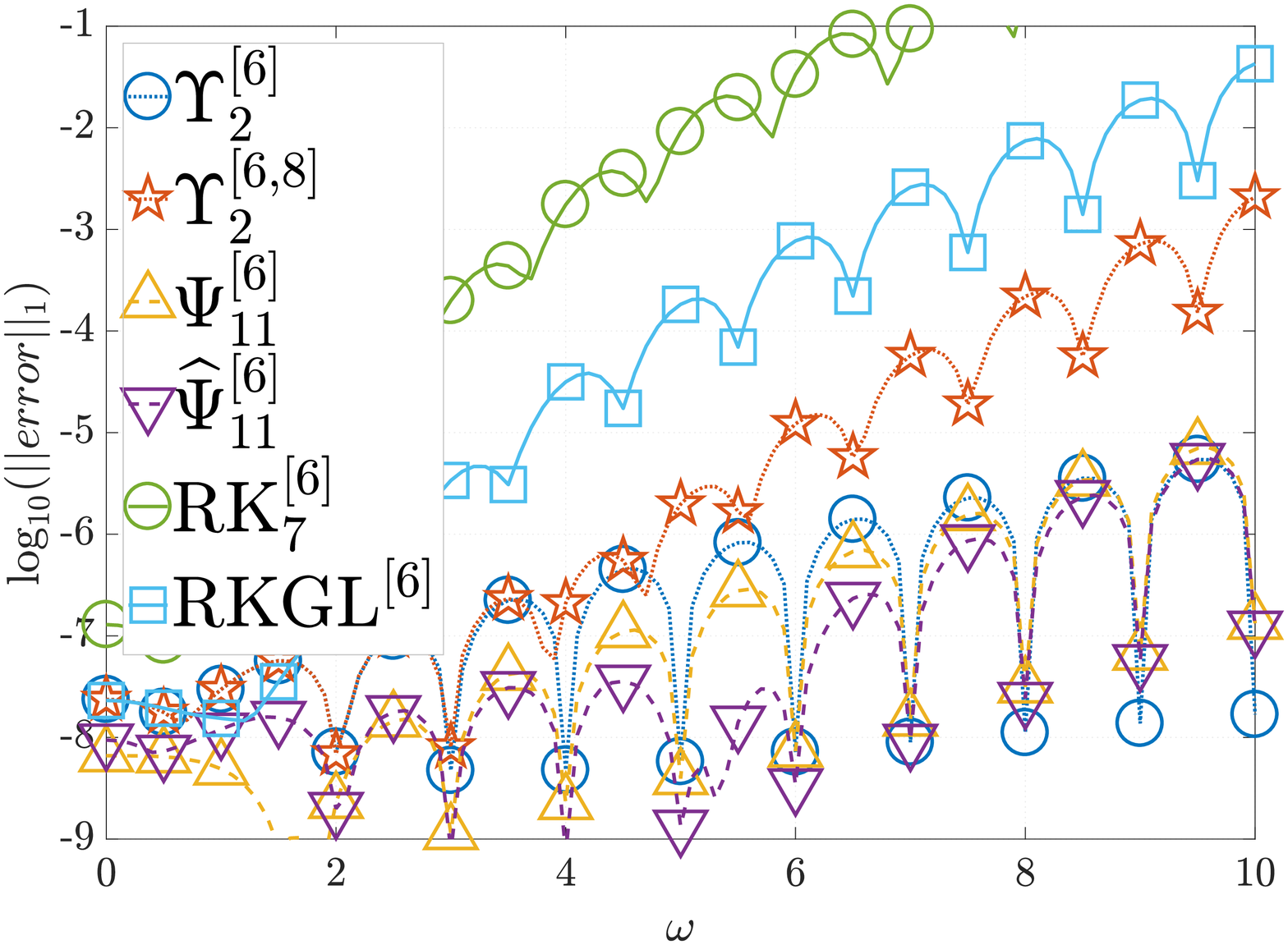}
			\caption{6th-order methods}
		\end{subfigure}			
		\caption{Error growth depending on $ \omega $ in the Mathieu \cref{ex:mathieu} on a logarithmic scale.}\label{fig:error-vs-omega}
	\end{figure}

	To illustrate how the~accuracy depends on the~frequency $\omega$, we take $\eps=1$, the~time step $ h=\pi/20 $ and~measure the~$ L_1 $-norm of the~error in the~fundamental matrix solution for $ \omega\in{[0,10]} $. The~results are shown in \Cref{fig:error-vs-omega}. We observe that the new methods show smaller error growth as frequency of the problem increases (oscillatory problems).
	%
	%
	
	%
	%
	%
	%
	
	\newpage
	\subsection{Hill equation}
	The second benchmark to~consider is the~matrix Hill equation:
	\begin{equation}\label{ex:hill}
		x''(t)+(A+B_1\cos2t+B_2\cos4t)x(t)=0
	\end{equation}
	where $ A,B_1,B_2 \in \mathbb{R}^{r\times r}$. We assume $ A=r^2 I +D $, where $ D $ is a~Pascal matrix:
	\[ 
		D_{1i}=D_{i1}=1,\: D_{ij}=D_{i-1,j}+D_{i,j-1}, 	\qquad 1<i,j\leq r. 
	\]
	We set $ B_1 =\eps I, B_2=\frac{1}{10}\eps I$, $ \eps=r $ and~$ \eps=\frac{1}{10} r $ and~compute solutions for $ r=5 $ and~$ r=7 $ on the~interval $ t\in[0,\pi] $, and then we measure the~$ L_1 $-norm of the~error in the~fundamental matrix solution at the final time.
	
\begin{figure}[!h]
		\begin{subfigure}[b]{0.5\textwidth} 	
			\includegraphics[width=1\linewidth]{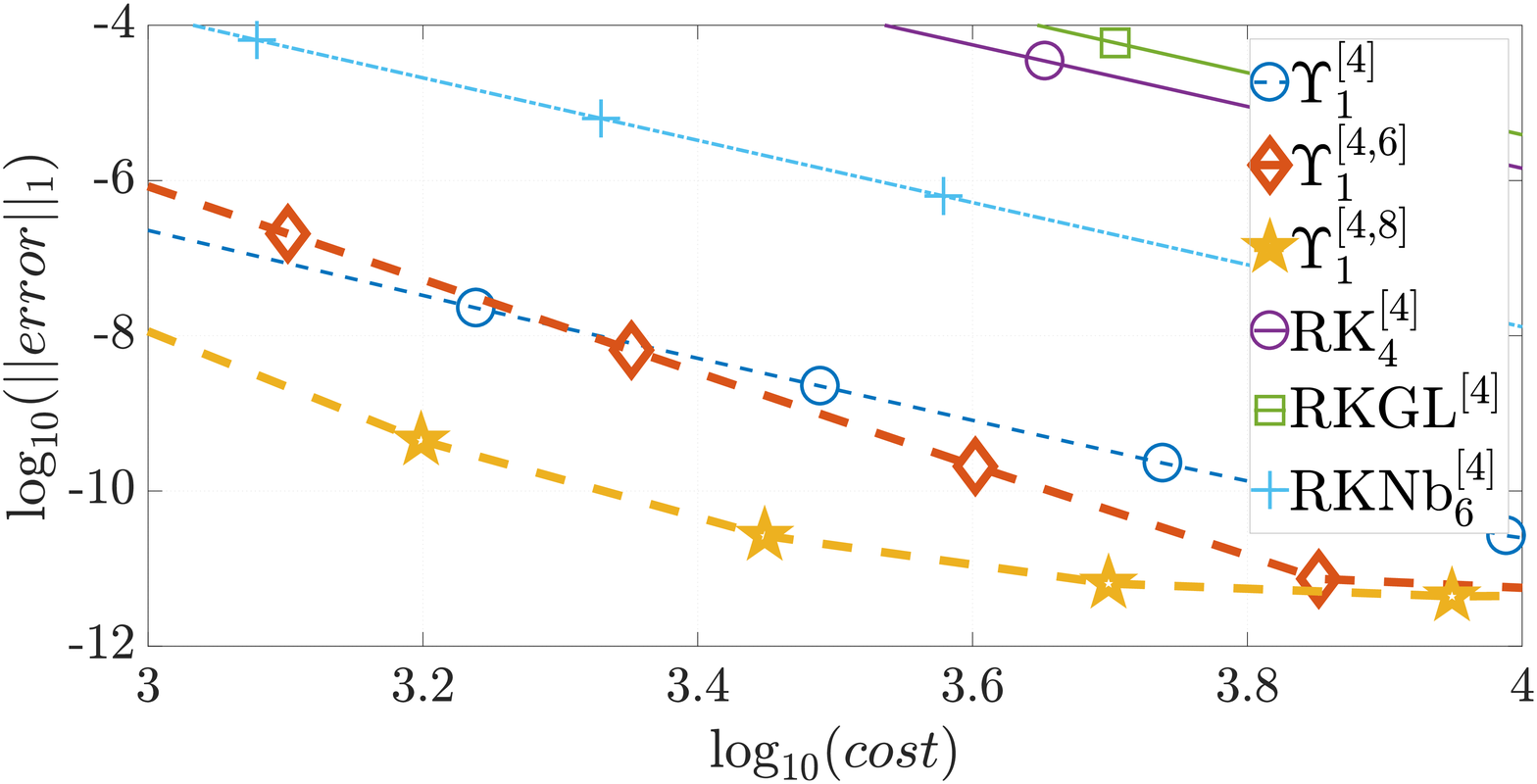}
			\caption{$ r=5, \eps=5/10 $}
		\end{subfigure}
		\begin{subfigure}[b]{0.5\textwidth} 	
			\includegraphics[width=1\linewidth]{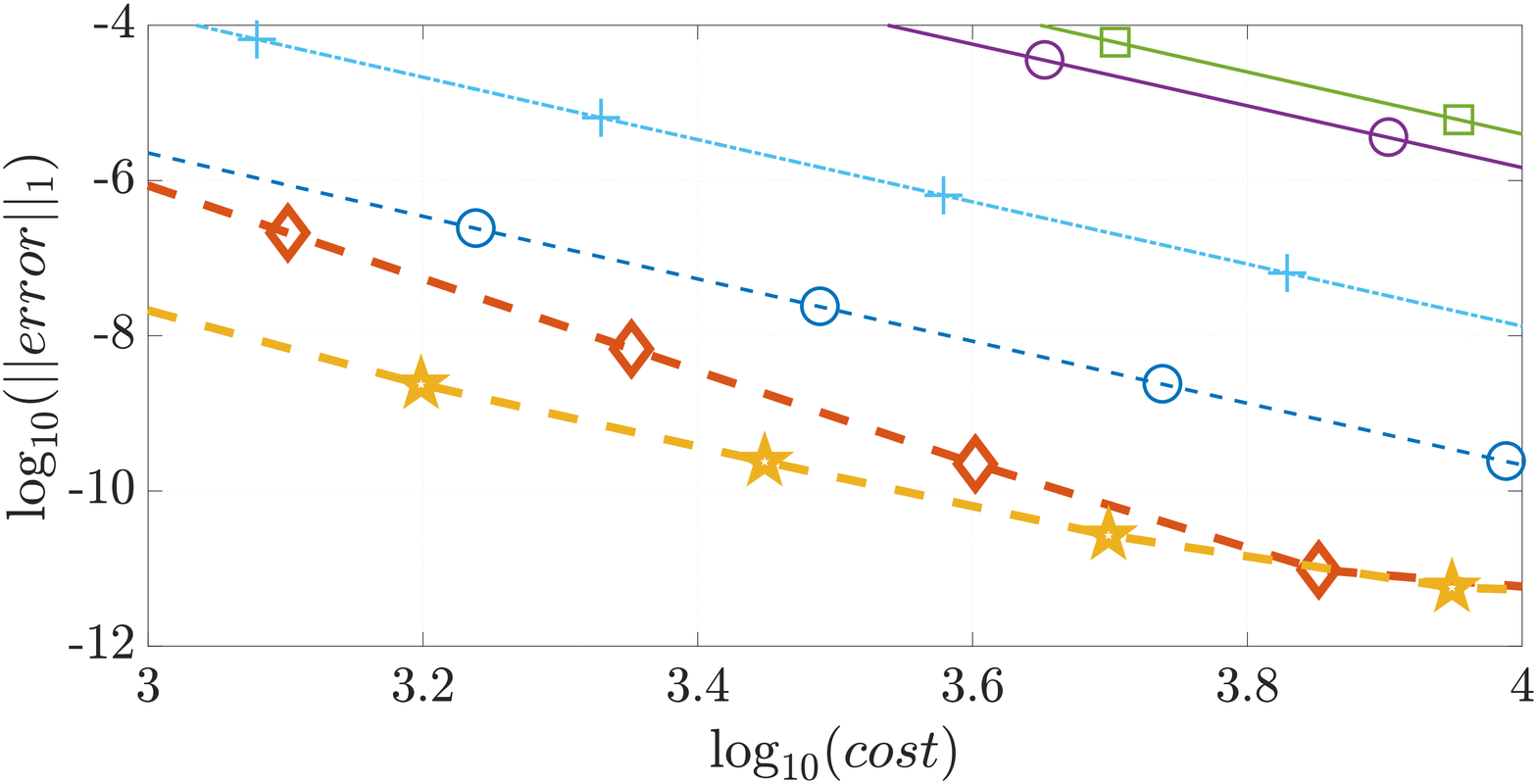}
			\caption{$ r=5, \eps=5 $}
		\end{subfigure}		
		\begin{subfigure}[b]{0.5\textwidth} 	
			\includegraphics[width=1\linewidth]{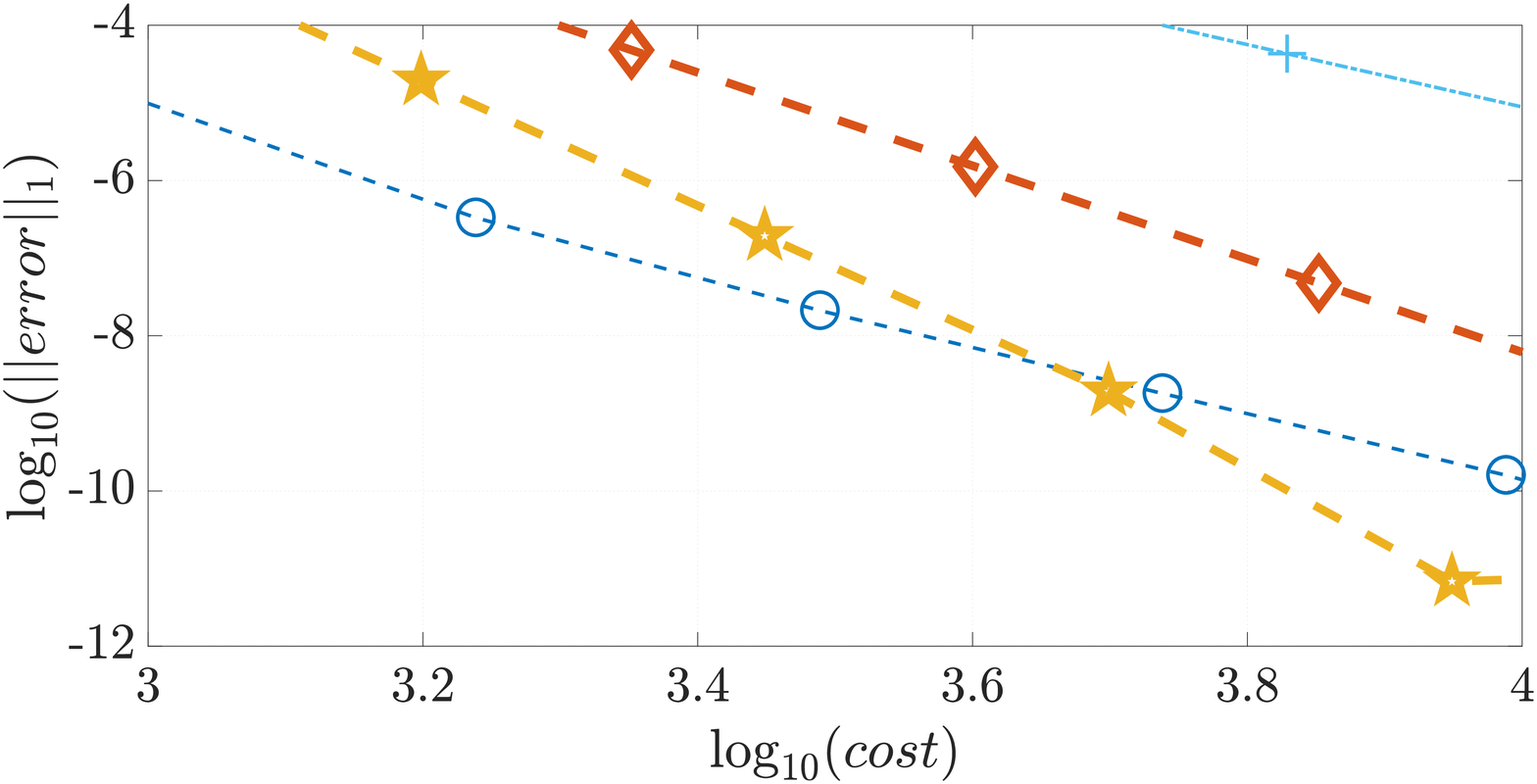}
			\caption{$ r=7, \eps=7/10 $}
		\end{subfigure}
		\begin{subfigure}[b]{0.5\textwidth} 	
			\includegraphics[width=1\linewidth]{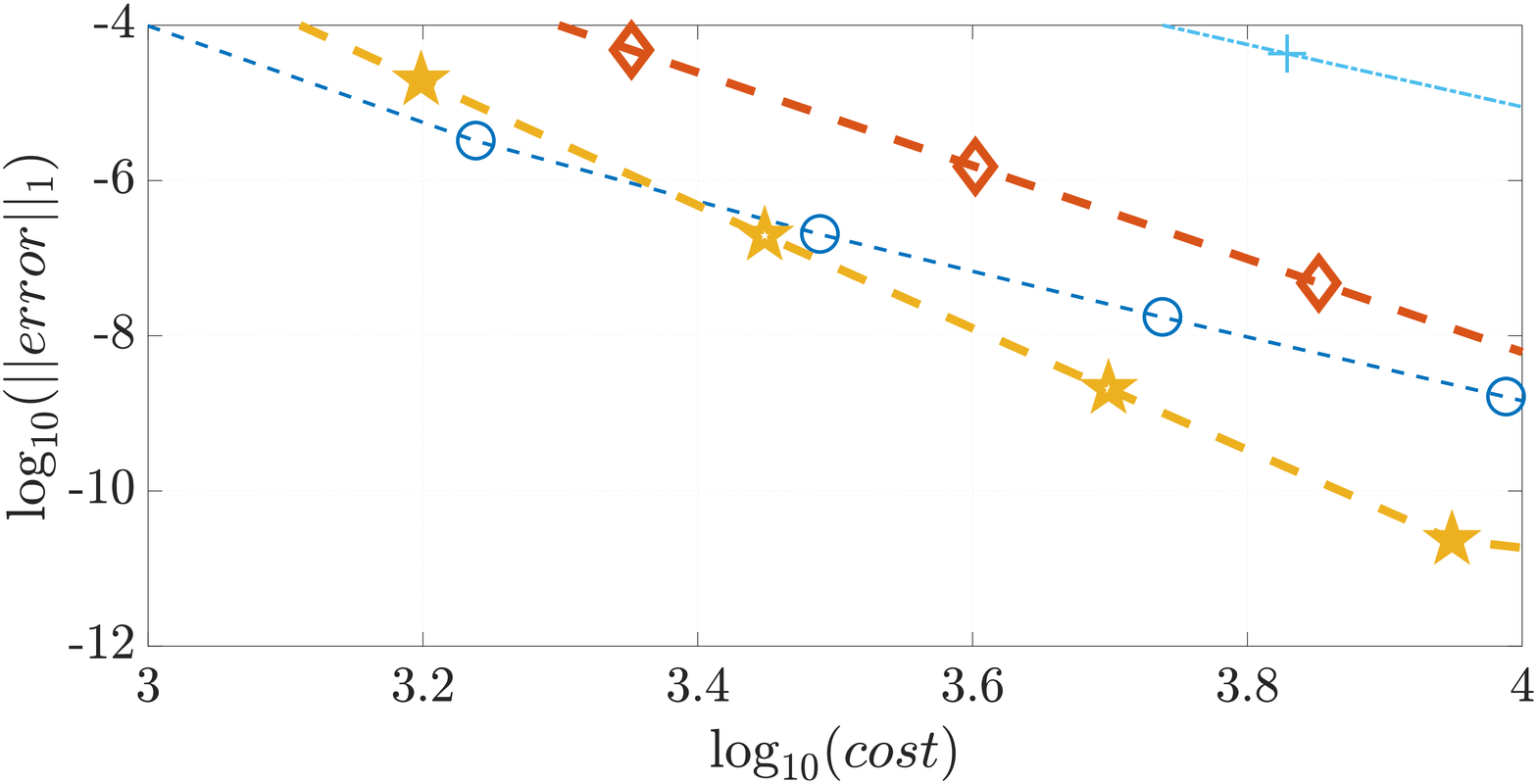}
			\caption{$ r=7, \eps=7 $}
		\end{subfigure}		
		\caption{The performance of the 4th-order methods for the Hill \cref{ex:hill} on a $ \log $--$ \log $ scale; $ cost = \mathcal{C}\times steps $.}
		\label{fig:hill_4th}
	\end{figure}

\begin{figure}[!h]
		\begin{subfigure}{0.5\textwidth} 	
			\includegraphics[width=1\linewidth]{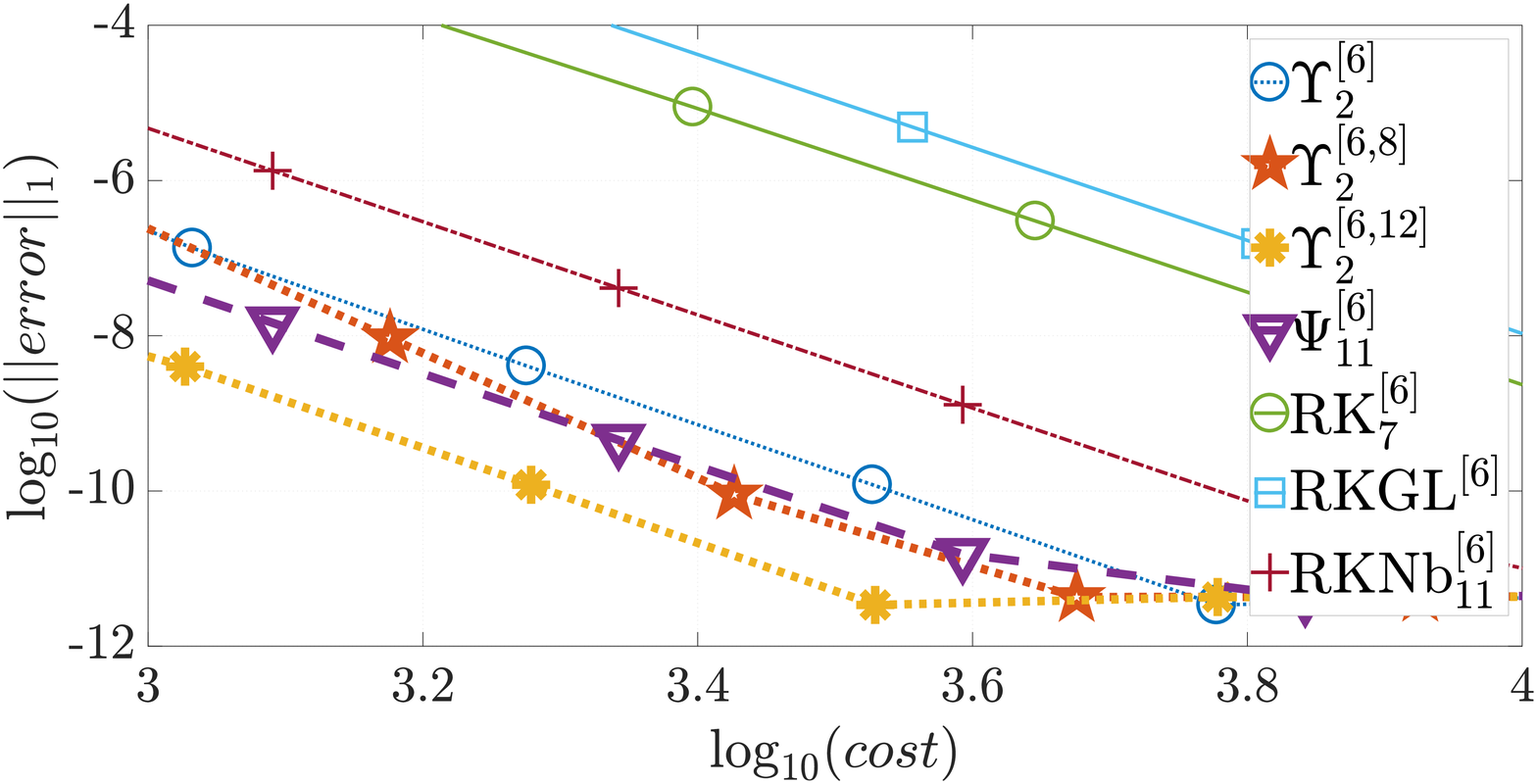}
			\caption{$ r=5, \eps=5/10 $}
		\end{subfigure}
		\begin{subfigure}{0.5\textwidth} 	
			\includegraphics[width=1\linewidth]{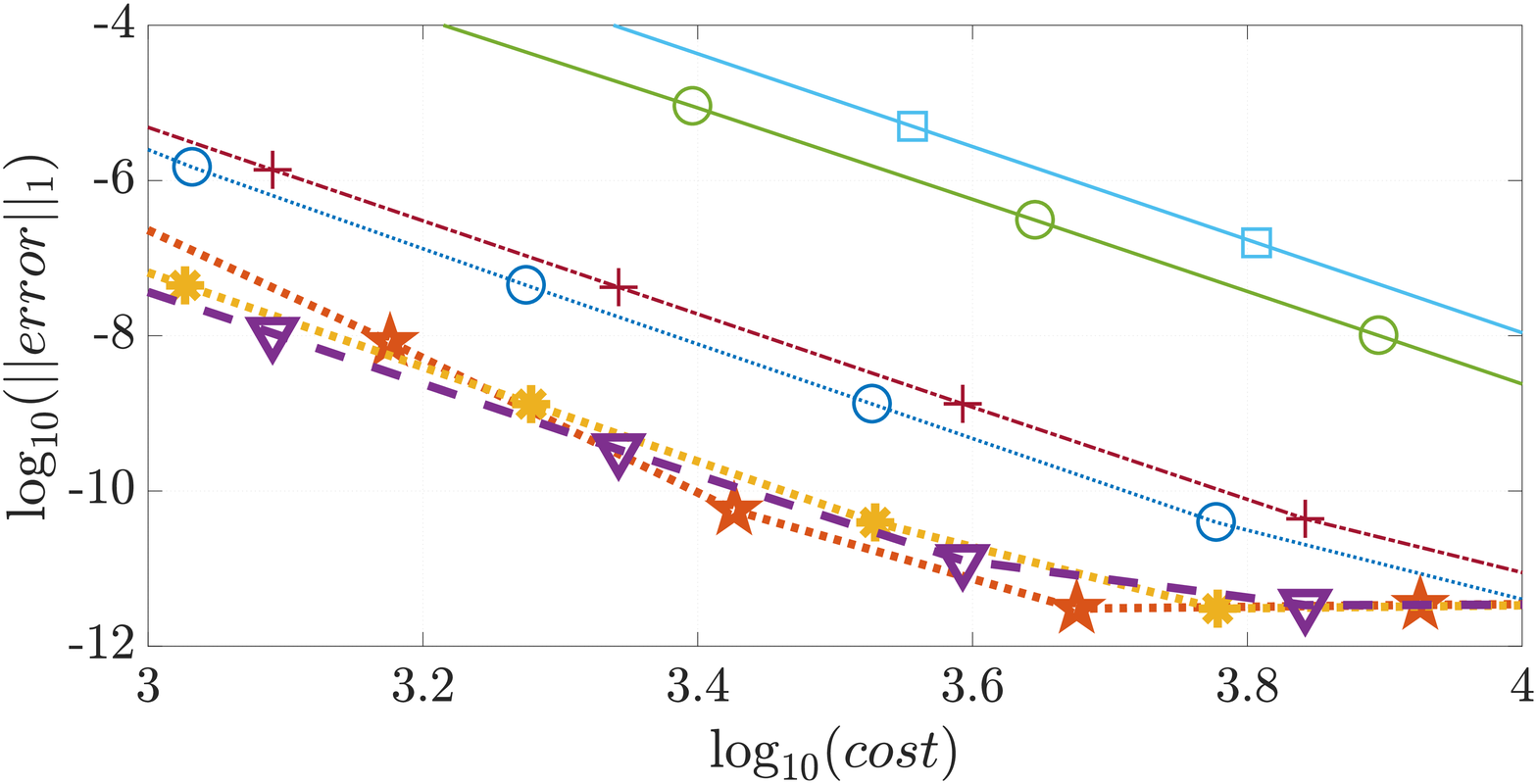}
			\caption{$ r=5, \eps=5 $}
		\end{subfigure}
		
		\begin{subfigure}[b]{0.5\textwidth} 	
			\includegraphics[width=1\linewidth]{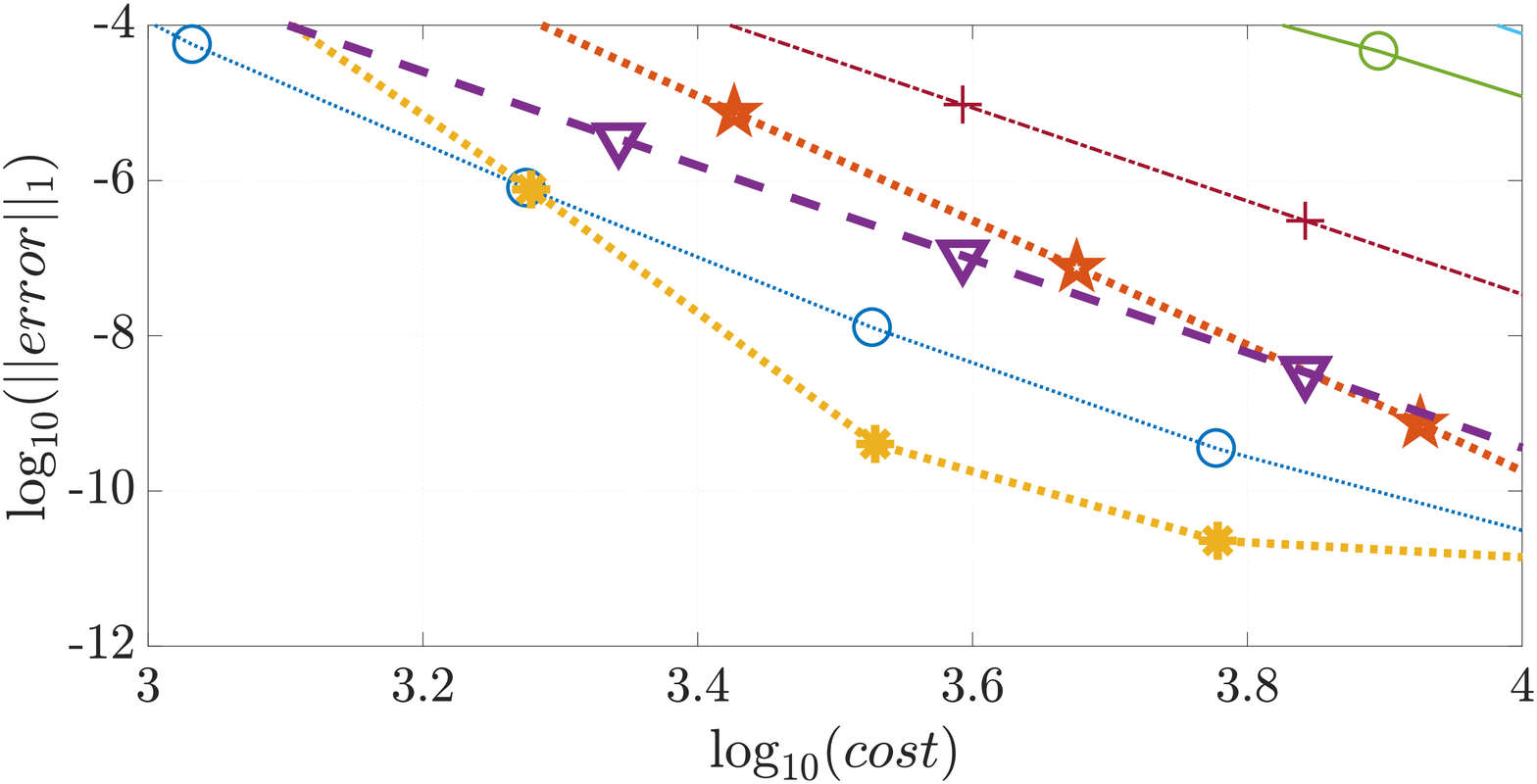}
			\caption{$ r=7, \eps=7/10 $}
		\end{subfigure}
		\begin{subfigure}[b]{0.5\textwidth} 	
			\includegraphics[width=1\linewidth]{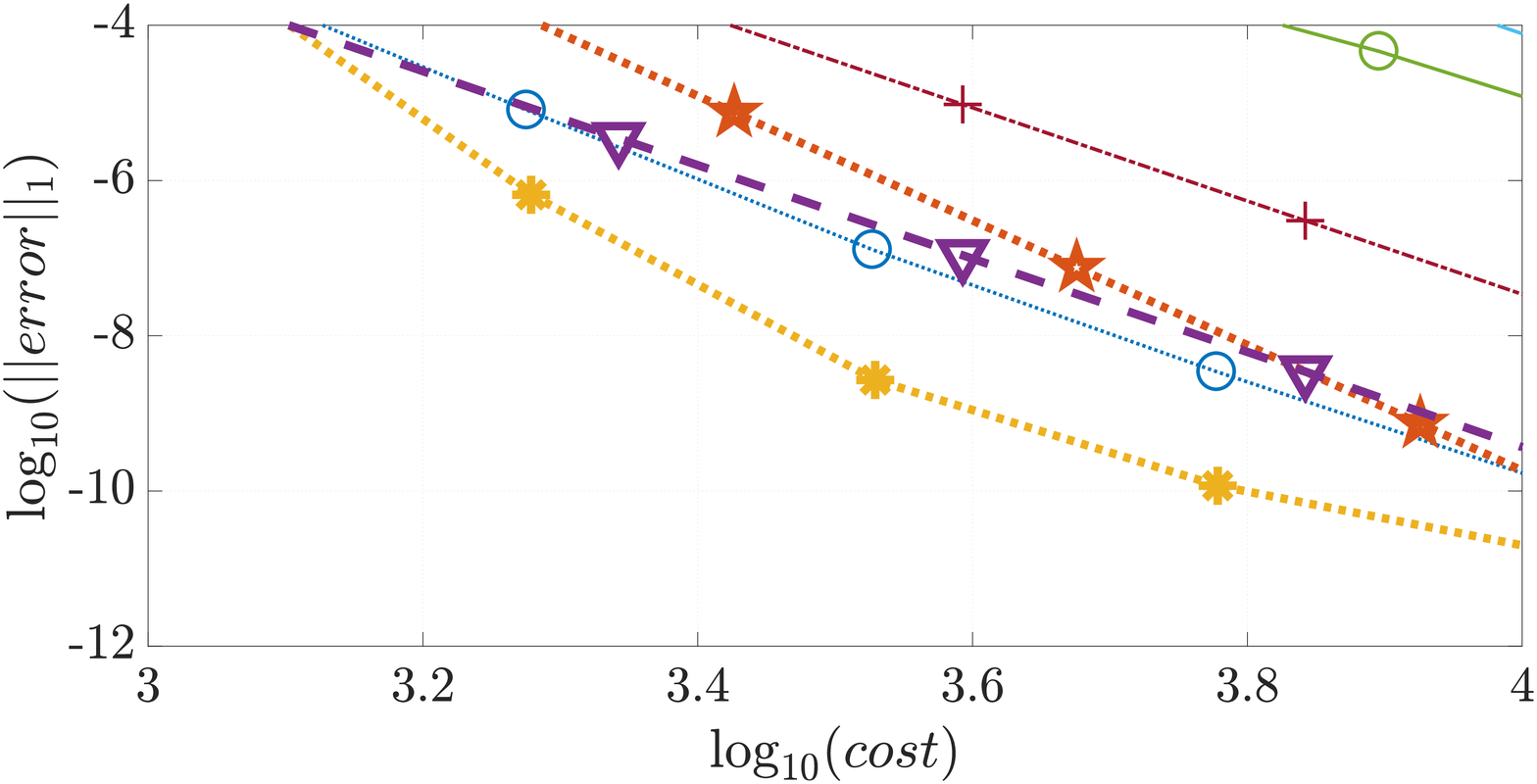}
			\caption{$ r=7, \eps=7 $}
		\end{subfigure}
		
		\caption{The performance of the 6th-order methods for the Mathieu \cref{ex:hill} on a $ \log $--$ \log $ scale; $ cost = \mathcal{C}\times steps $.}
		\label{fig:hill_6th}
	\end{figure}

	In \cite{bader16sif}, it was determined that for matrix Hill-type problems $\Upsilon_{k}^{[p]} $ perform no worse than $ \mathrm{RKNb}_{11}^{[6]} $. \Cref{fig:hill_4th,fig:hill_6th} demonstrate that, thanks to~lower computational cost, the new $\Upsilon_{1}^{[4,q]}, \Upsilon_{2}^{[6,q]} $ and~$ \Psi_{11}^{[6]} $ methods consistently produce better results in both oscillatory (larger $ r $) and~nearly autonomous (small $ \eps $) cases.
	%
	%
	%
	%
		\subsection{Time-dependent wave equation}
	
	To analyse the~performance of the~methods that only involve matrix--vector products we consider the~following trapped wave equation
	\begin{equation}\label{eq:wave_eq}
		\partial_t^2 u = \partial_x^2 u - \big( x^2+g(x,t) \big) u,
		\quad x\in\mathbb{R},\ t\geq 0,
	\end{equation}
	equipped with initial conditions $u(x,0)=\sigma \e^{-x^2/2}$, and~$u_t(x,0)=0$. The~solution for $g(x,t)=0$ can be easily be obtained by separation of variables and it is given by $u_0(x,t)=\sigma \cos(t)\e^{-x^2/2}$.
	
	When an external interaction appears, $g(x,t)\neq 0$, the~equation has no analytical solution in general and one has to~consider, for example, a~numerical scheme.
	We assume the~solution is confined to a~region $x\in[x_0,x_N]$ and hence the~solution and all spatial derivatives vanish at these boundaries.
	This allows us to treat the problem as periodic and spectral methods can be used. 
	We divide the~spatial region into $N$ intervals of length $\Delta x=(x_N-x_0)/N$ and, after spatial discretization, we obtain an equation similar to~\cref{hill} that we write as the~first order system \cref{firstordersystem} where $z=(v,w)^T$ and $v_i(t)\approx u(x_i,t), w_i(t)\approx u_t(x_i,t)$.
	
\begin{figure}[!h]	
		\begin{subfigure}[b]{0.5\textwidth}	
			\centering
			\includegraphics[width=1\linewidth]{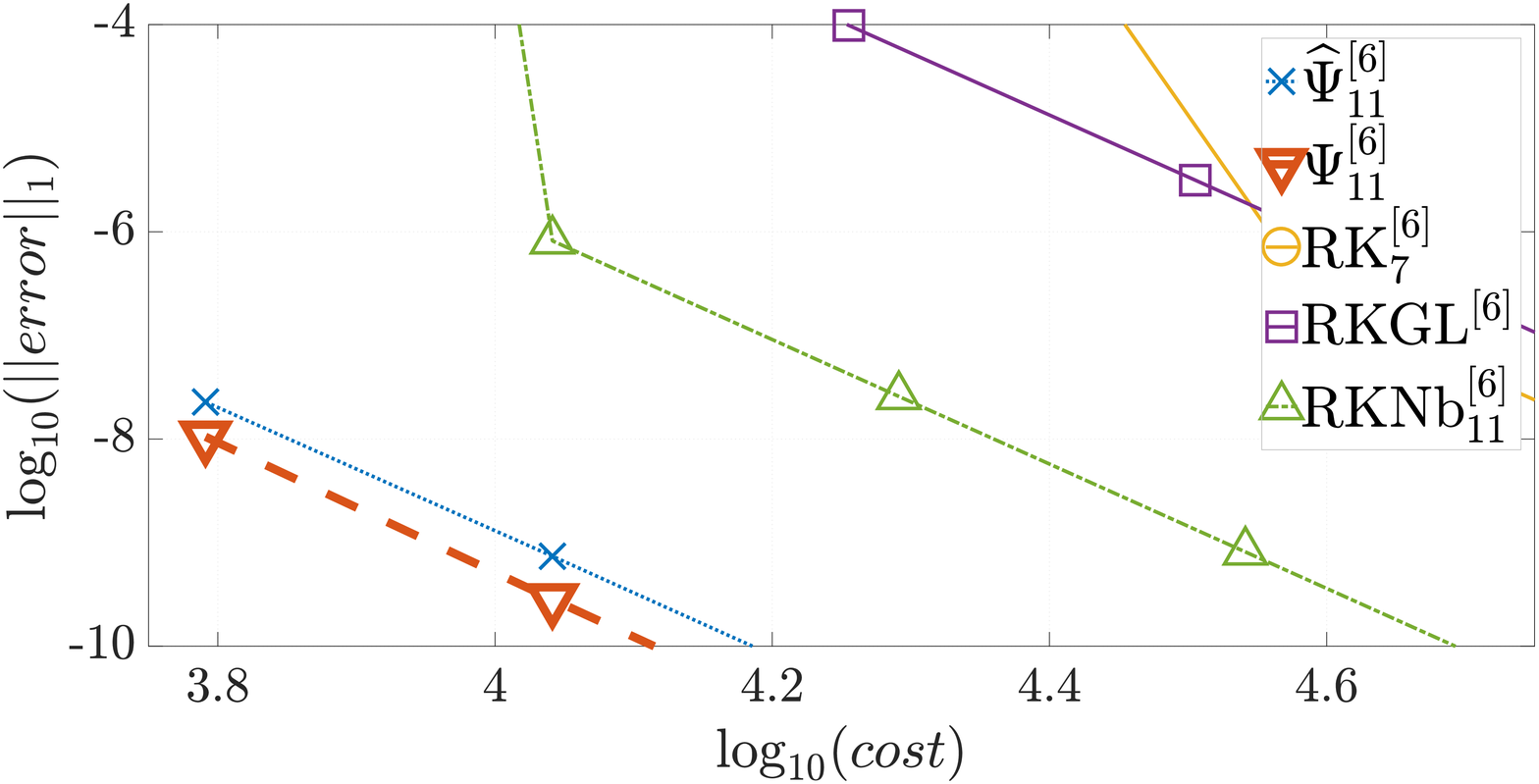}
			\caption{$ \delta=0.2,\ \eps=0.2 $}
		\end{subfigure}
		\begin{subfigure}[b]{0.5\textwidth}	
			\centering
			\includegraphics[width=1\linewidth]{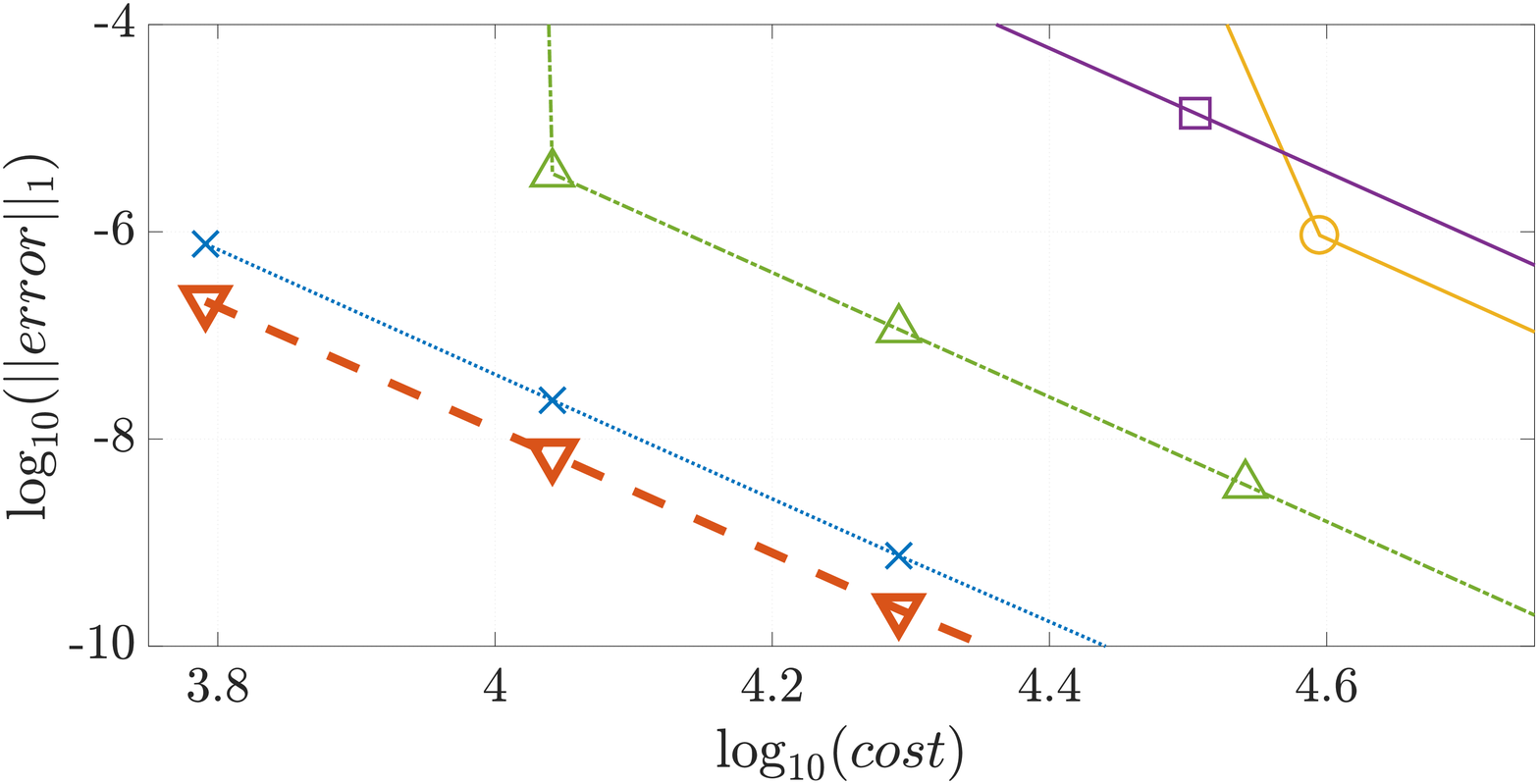}
			\caption{$ \delta=0.2,\ \eps=0.4 $}
		\end{subfigure}
	\begin{subfigure}[b]{0.5\textwidth}	
		\centering
		\includegraphics[width=1\linewidth]{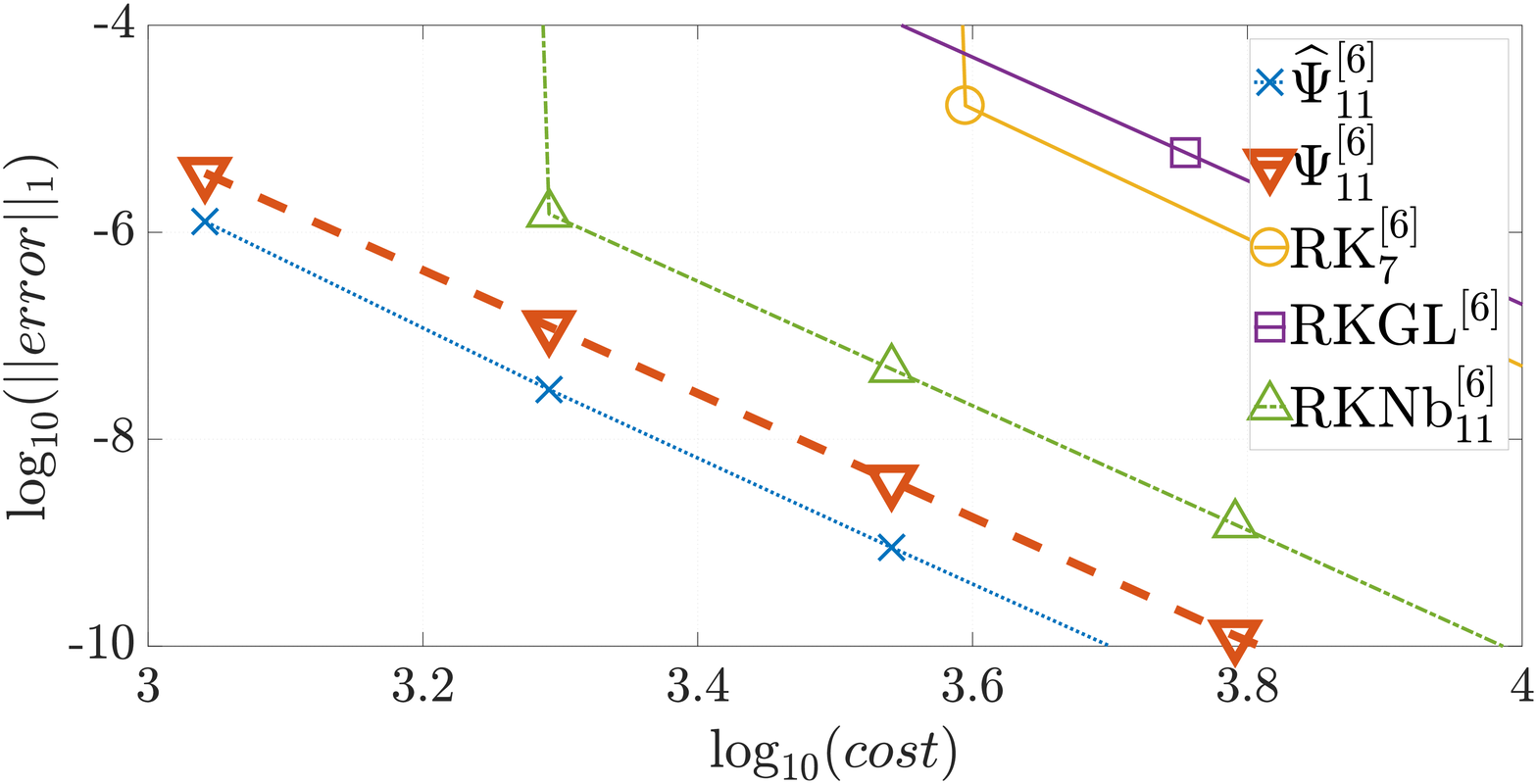}
		\caption{$ \delta=1,\ \eps=0.1 $}
	\end{subfigure}
	\begin{subfigure}[b]{0.5\textwidth}	
		\centering
		\includegraphics[width=1\linewidth]{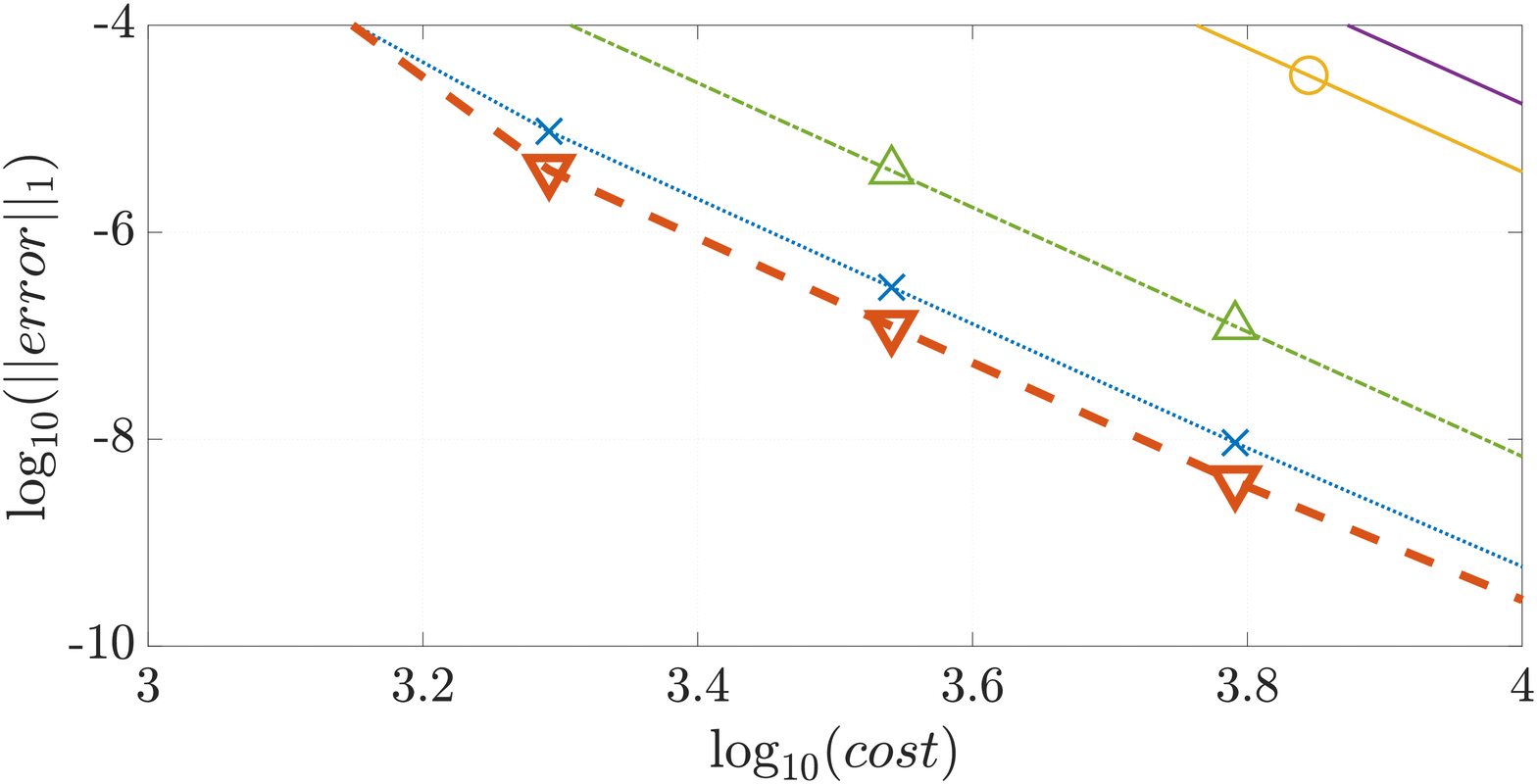}
		\caption{$ \delta=1,\ \eps=0.5 $}
	\end{subfigure}
		\caption{The performance of the 6th-order methods for the wave \cref{eq:wave_eq} on~a~$ \log $--$ \log $ scale; $ cost = \mathcal{C}\times steps $.}
		\label{fig:wave_eq_6th}
	\end{figure}

	For the~numerical experiments we take $N=128, x_0=-10, x_{N}=10$ and the~external interaction
	\[
		g(x,t)= \eps\cos (\delta\, t) x^2 .
	\]
	We take $\delta\in \{0.2,\:1\}$, $\eps\in\{0.1,0.2,\:0.4,\:0.5\}$, and integrate for the period $t\in[0,20\pi/\delta]$. The reference solution is obtained numerically with a~sufficiently small time step and we measure the~$L_1$-norm of the~solution versus the~number of matrix\,--\,vector products for each method. The~results are shown in the \Cref{fig:wave_eq_6th} where $ \Psi_{11}^{[6]}$, with only three evaluations of the time-dependent functions per step, shows the~best performance for smooth time-dependencies and it is nearly the best one when such time-dependency increases. 
	%
	%

	\section{Conclusions}
	
	        Starting from the well know Magnus expansion for linear time-dependent differential equations, in this work we have presented two families of time-integrators especially designed for the second order linear system (\ref{hill}).         
        Both
share with the exact solution relevant qualitative properties (in particular, its symplectic character) and thus provide an accurate description
of the system over long time intervals, due to their favourable error propagation, but they differ in significant aspects. Whereas Magnus-decomposition
methods are addressed to problems with small to moderate dimensions where the numerical computation of matrix--matrix products (and thus
the fundamental matrix) is feasible, Magnus--splitting methods are advantageous when the dimension of the system is exceedingly high, as is the
case when (\ref{hill}) results from a linear time-dependent wave equation previously discretized in space. The numerical experiments reported
here clearly illustrate this difference: whereas the first family performs more efficiently on the scalar Mathieu equation and a low-dimensional
and oscillatory matrix Hill equation, it is the second class of integrators which shows a better behaviour on a discretised trapped wave equation,
where algorithms based only in matrix-vector products are advisable. 

Magnus-decomposition methods can be considered as an adaptation of the schemes presented in \cite{bader16sif} involving appropriate 
truncated expansions of the otherwise computationally costly matrix exponentials initially present in them and useful combinations of matrices. 
Magnus--splitting methods, on the other hand, belong to the class of integrators analysed in \cite{blanes07smf,blanes12smi}: taking as a starting
point an efficient splitting method for the autonomous case, the scheme is formed by averaging the matrix $M(t)$ at each step (with different
weights). The new coefficients are obtained by solving the additional order conditions arising from the time-dependency. Compared with
the general case, the number of order conditions is considerably reduced, which allows to get new schemes within this family in a rather
straightforward way.

	\section*{Acknowledgements}
	
	Bader, Blanes, Casas and Kopylov acknowledge the~Ministerio de Econom\'ia y Competitividad (Spain) for financial support through the~coordinated project MTM2013-46553-C3-3-P. Additionally, Kopylov has been partly supported by fellowship GRISOLIA/2015/A/137 from the~Generalitat Valenciana.


\section*{References}
	\begin{thebibliography}{99}
		%
		\bibitem{almohy15naf} A.H. Al-Mohy, N.J. Higham and S.D. Relton,
		\newblock New Algorithms for Computing the~Matrix Sine and Cosine Separately or Simultaneously,
		SIAM J. Sci. Comput. {\bf 37} (1) (2015), A456-A487.
		%
		\bibitem{alonso16eaa} P. Alonso, J. Ib\'an\~ez, J. Sastre, J. Peinado, and E. Defez, 
		Efficient and accurate algorithms for computing matrix trigonometric functions,
		\textit{J. Comp. Appl. Math.} \textbf{309}, pp. 325--332, 2017.
		%
		\bibitem{bader14stp}		P. Bader, S. Blanes, Solving the perturbed quantum harmonic oscillator in imaginary time using splitting methods with complex coefficients, in F. Casas, V. Martinez (eds.), Advances in Differential Equations and Applications, SEMA SIMAI Springer Series 4, (2014), pp. 217--227.
		
		\bibitem{bader16eni} P. Bader, S. Blanes, F. Casas, and~E. Ponsoda,
		Efficient numerical integration of Nth-order non-autonomous linear differential equations.
		\textit{J. Comp. Appl. Math.} \textbf{291}, pp. 380--390, 2016.
		%
		\bibitem{bader16sif} P. Bader, S. Blanes, E. Ponsoda and~M. Seydao\u{g}lu,
		Symplectic integrators for the~matrix Hill equation.
		\textit{J. Comp. Appl. Math.} In press.
		%
		\bibitem{bader16emf} 
		P. Bader, A. Iserles, K. Kropielnicka, and P. Singh,
		Efficient methods for linear Schr\"odinger equation in the
		semiclassical regime with time-dependent potential, 
		Proc. R. Soc. a~\textbf{472} (2016) 20150733.
		%
		\bibitem{Bernstein2009}	D.~S. Bernstein.
		\newblock Matrix mathematics: theory, facts, and~formulas, vol.~1, 2009.
		%
		\bibitem{blanes16aci} S. Blanes, F. Casas, 
		A Concise Introduction to~Geometric Numerical Integration. 
		CRC Press, Boca Raton, (2016).
		%
		\bibitem{blanes07smf} S. Blanes, F. Casas, and~A. Murua, 
		Splitting methods for	non-autonomous linear systems, 
		\textit{Int. J. Comput. Math.}, \textbf{84} (2007), pp. 713-727.
		%
		\bibitem {blanes12smi} S. Blanes, F. Casas, and~A. Murua, 
		Splitting methods in the~numerical integration of non-autonomous dynamical systems, \textit{RACSAM}, \textbf{106} (2012), pp. 49-66.
		%
		\bibitem{blanes09tme}	S.~Blanes, F.~Casas, J.~A.~Oteo, J.~Ros.
		\newblock the~Magnus expansion and~some of its applications.
		\newblock {Physics Reports}, \textbf{470}, (2009), pp. 151--238.
		%
		\bibitem{blanes02psp}	S. Blanes and~P.C. Moan, 
		Practical Symplectic Partitioned Runge-Kutta and~Runge-Kutta-Nyström Methods, 
		J. Comput. Appl. Math., \textbf{142} (2002), pp. 313-330.
		%
		\bibitem{dragt15lmf} A. J. Dragt,
		Lie Methods for Nonlinear Dynamics with Applications to~Accelerator Physics,
		University of Maryland, 2015.
		%
		\bibitem{drewsen00hlp} M. Drewsen and A. Br{\o}ner,
		Harmonic linear Paul trap: Stability diagram and effective potentials,
		{Phys. Rev. A} {62} (2000) 045401.
		%
		\bibitem{hairer06gni}{E. Hairer, Ch. Lubich, and G. Wanner}, 
		{Geometric Numerical Integration}, 2nd ed., Springer, Berlin, 2006.
		%
		\bibitem{iserles99ots} A. Iserles and S.P. N{\o }rsett, 
		On the~solution of linear differential equations in Lie groups, 
		{Philos. Trans. Royal Soc. London Ser. A} {357} (1999), pp. 983--1019.
		%
		\bibitem{leimkuhler05shd} {B. Leimkuhler and S. Reich}, 
		{Simulating Hamiltonian Dynamics,} Cambridge University Press, Cambridge, 2005.
		%
		\bibitem{magnus54ote}	W.~Magnus.
		\newblock On the~exponential solution of differential equations for a~linear operator.
		\newblock {Comm. Pure and~Appl. Math.}, \textbf{VII}, (1954) pp. 649--673.
		%
		\bibitem{magnus66he} W. Magnus and~S.Winkler, 
		Hill equation. Wiley, New York, 1966. 
		%
		\bibitem{major05cpt} F. G. Major, V. N. Gheorghe, and G. Werth, 
		Charged Particle Traps. Physics and Techniques of Charged Particle Field Confinement, Springer, 2005.
		%
		\bibitem{moan98eao} P.C. Moan,
		Efficient Approximation of {S}turm--{L}iouville Problems Using {L}ie-Group Methods. DAMTP, Tech. Report 1998/NA11, University of Cambridge, United Kingdom (1998).
		%
		\bibitem{mclachlan65taa} N.W. McLachlan, 
		Theory and application of Mathieu functions, Dover, New York, 1964.
		%
		\bibitem{munthekaas99cia} H. Munthe-Kaas and B. Owren, 
		Computations in a~free Lie algebra, Philos. Trans.
		Royal Soc. London Ser. A, 357 (1999), pp. 957–981.
		%
		\bibitem{paul90etf} W. Paul, 
		Electromagnetic traps for charged and neutral particles,
		{Rev. Modern Phys.} {62} (1990) 531--540.
		%
		\bibitem{sanzserna94nhp} {J. M. Sanz-Serna and M. P. Calvo}, 
		{Numerical Hamiltonian Problems,} Chapman and Hall, London, 1994.
		%
		\bibitem{turner98fpr} K. L. Turner, S. A. Miller, P. G. Hartwell, N. C. MacDonald, S. H. Strogatz, and S. G. Adams,
		Five parametric resonances in a~microelectromechanical system,
		Nature {396} (1998) 149--152.
		%
		\bibitem{zanna99car} A. Zanna, 
		Collocation and relaxed collocation for the~Fer and the~Magnus expansions,
		SIAM J. Numer. Anal., 36 (1999), pp. 1145--1182.
		
		
		
	\end{thebibliography}

\end{document}